\documentclass[11pt]{article}
\usepackage{amsmath, amsthm, amssymb, setspace, mathrsfs, mathpazo, enumitem, mathtools, xypic}
\makeatletter
\renewcommand\section{\@startsection{section}{1}{\z@}%
 						{-3.5ex \@plus -1ex \@minus -.2ex}
						{2ex \@plus.2ex}
						{\large\bfseries}}
\renewcommand\subsection{\@ifstar
						{\setcounter{subsection}{\value{equation}}
					\@startsection{subsection}{2}{\z@}
                          {1.75ex \@plus.5ex \@minus.2ex}%
                           {-.4em}		
					\textit*}
					{\setcounter{subsection}{\value{equation}}
						\stepcounter{equation}
					\@startsection{subsection}{2}{\z@}
                          {1.75ex \@plus.5ex \@minus.2ex}%
                           {-.4em}		
					\textit}}
\def\@seccntformat#1{\@ifundefined{#1@cntformat}%
	{\csname the#1\endcsname\quad} 
	{\csname #1@cntformat\endcsname}} 
\def\section@cntformat{\thesection.~} 
\def\subsection@cntformat{(\thesubsection)\ }
\renewcommand*\l@section{\mdseries\small\@dottedtocline{1}{1.5em}{2em}}

\newcommand{\xdbheadrightarrow}[2][]{%
  \ext@arrow 0099\xdbheadfill@{#1}{#2}}%
\newcommand{\xdbheadfill@}{%
  \arrowfill@\relbar\relbar{\mathrel{\vphantom{\rightarrow}\smash{\twoheadrightarrow}}}}
\makeatother

\numberwithin{equation}{section}
\theoremstyle{plain}
\newtheorem{maintheorem}{Theorem}
\swapnumbers
\newtheorem{theorem}[equation]{Theorem}
\newtheorem*{theorem*}{Theorem}
\newtheorem{corollary}[equation]{Corollary}
\newtheorem{lemma}[equation]{Lemma}
\newtheorem{proposition}[equation]{Proposition}

\theoremstyle{definition}
\newtheorem{definition}[equation]{Definition}
\theoremstyle{remark}
\newtheorem{remark}[equation]{Remark}

\newcommand{\script}{\mathscr}
\newcommand{\cC}{\script{C}}

\newcommand{\caH}{\script{H}}
\newcommand{\cK}{\script{K}}
\newcommand{\caR}{\script{R}}

\newcommand{\cW}{\script{W}}

\newcommand{\frg}{\mathfrak{g}}
\newcommand{\frh}{\mathfrak{h}}

\newcommand{\bC}{\mathbb{C}}
\newcommand{\bH}{\mathbb{H}}

\newcommand{\bP}{\mathbb{P}}
\newcommand{\bR}{\mathbb{R}}
\newcommand{\bQ}{\mathbb{Q}}
\newcommand{\bS}{\mathbb{S}}

\newcommand{\bZ}{\mathbb{Z}}
\newcommand{\vep}{\varepsilon}
\newcommand{\Ind}{\mathrm{Ind}}
\newcommand{\pz}{(\!(z)\!)}
\newcommand{\bz}{[\![z]\!]}

\newcommand{\bmu}{\boldsymbol{\mu}}
										
\begin{document}
\title{\textbf{Coulomb branches for quaternionic representations}}
\author{Constantin Teleman}
\date{\today}
\maketitle
\begin{quote}
\abstract{
\noindent 
I describe the \emph{Chiral rings} $\caR_{3,4}$ for $3$D, $N=4$ supersymmetric $G$-gauge theory and 
matter fields in quaternionic representations~$E$: first, by incorporating \emph{twisted} real 
structures in the construction of~\cite{bfn}, and second, more explicitly, by Weyl group descent 
from the maximal torus. A topological obstruction  
$w_4(E)$ (modulo squares) appears for $\caR_3$; 
a secondary obstruction $\eta\cdot E$, appears for $\caR_4$. The two combine to a gauging 
obstruction of~$E$ by~$G$ in $4$D, $N=2$ supersymmetry, enhancing Witten's obstruction. I classify the obstructions 
for connected~$G$. Freedom of 
the chiral rings over the Toda bases reduces calculations to the case of semi-simple rank~$1$. For 
representations whose weights include the roots of~$G$, an Abelianization formula describes 
the $\caR$ in terms of the maximal torus and the Weyl group. My approach an alternative and generalization 
to a recent paper \cite{betal}.}
\end{quote}


\section*{}
\begin{center}
\begin{minipage}[t]{12cm}
\tableofcontents
\end{minipage}
\end{center}
 
\section*{Introduction}

For a compact Lie group $G$ and a quaternionic representation $E$ 
there are expected to be (singular) 
hyper-K\"ahler \emph{Coulomb branches} of the moduli of vacua for ($3$D, $N=4$ supersymmetric) gauge 
theory, with matter fields in $E$. Arising, as they do, from four dimensions by dimensional reduction 
along a line or a circle, they come in two versions, $\cC_3(G;E), \cC_4(G;E)$, respectively. 
In a familiar pattern, they are associated to ordinary cohomology and to $K$-theory: the equivariant 
parameter for circle rotation becomes the (inverse of the) Bott periodicity generator. 

The complex varieties $\cC_{3,4}(G;0)$ are the total spaces for the Toda integrable system and its finite 
difference version, respectively, with bases $\mathrm{Spec}\, H^*_G, K^*_G$ of a point. They were thoroughly 
studied in \cite{bfm}, although not, at the time, related to $3$D gauge theory. That connection was spelt out later 
\cite{telicm}, but special cases were known in the physics literature, such as $\mathrm{SU}(2)$ in \cite{sw}. 

For \emph{polarized} representations $E$ (symplectic doubles of complex ones), the first general construction 
of the Poisson algebras $\caR_{3,4}(G;E)$ of algebraic functions on the modified varieties $\cC_{3,4}(G;E)$ --- called \emph{chiral rings} in the physics literature\footnote{The term \emph{Coulomb branch} is sometimes 
understood to incorporate the hyperk\"ahler metric. That is known for $\cC_3(G;0)$, via the Nahm equations, 
using results of Bielawski \cite{bi}, and for abelian $G$ \cite{bdg}. A general construction for $\cC_3$ was recently proposed in \cite{bif}.}--- was given by Braverman, Finkelberg and Nakajima \cite{bfn}; 
see also \cite{bdg} for a physics perspective. These varieties should control gauge theory ``with matter fields 
in $E$.'' While a precise formulation of that goal is not yet at hand, their subsequent treatment in \cite{telc} proposed the moral teaching (partly generalized in the present paper) that the \emph {Chiral rings with matter 
are formal consequences of the 2$D$ gauged linear Sigma model (GLSM)}, viewed as a boundary condition. 

Following that lead, and the proposed relation with quantum and symplectic cohomology, a conceptual interpretation 
of the GLSM construction was recently established by Gonzalez, Mak and Pomerleano in \cite{gmp}: 
the chiral ring $\caR_3(G;E)\subset \caR_3(G;0)$ is the subring  which preserves the equivariant quantum cohomology 
lattice of a polar half of $E$ within its symplectic cohomology. This may well be the first substantive 
theorem --- as opposed to construction, or definition --- relating Coulomb branches to pre-existing mathematics.

In this paper, I describe the Poisson rings $\caR_{3,4}(G;E)$  believed to underlie the $\cC_{3,4}$ in two 
different ways. First, I adapt the construction of \cite{bfn} by exploiting a \emph{twisted form} of real 
structures, which untwists on ($K$-)homology following the removal of some obstructions. Second, I 
adapt the ``mirror B-model'' construction from \cite{telc}. The original version glues two copies 
of the Toda space along a vertical shift by the rational section $\exp(d\Psi)$, for the 
superpotential $\Psi$ of the GLSM. The analogue here is more complicated: $\Psi$ is 
now Weyl-multi-valued.    

The first construction originates topologically in an extension of the group $\Omega G$ by the multiplicative 
\emph{monoid} of ($K$-)homology, rather than by the group of units: the Euler class of the extension 
\begin{equation}\label{firstext}
BSO \rightarrowtail BU \twoheadrightarrow \Omega Sp^{w_4},
\end{equation} 
pulled back to $\Omega G$ by the Atiyah index map for $E$. (We obtain $Sp^{w_4}$ from~$Sp$ by 
killing $w_4$, the obstruction to constructing $\caR_3$.) A similar story applies to $\caR_4$, but 
with~$BSpin^c$ replacing~$BSO$. Ignoring the fiber of the extension and forcibly lifting to~$BU$ 
suspends $\Omega G$ by a (locally defined) linear space~$R_E$ with jumping fibers, playing the role of a 
\emph{half} of the index sheaf. The $J$-homomorphism converts it into a constructible coefficient 
system for ($K$-)homology. This $BU$-lift provides the `quantum square root' substituting for a missing polar 
half of~$E$. This can be extended to other generalized cohomology theories, subject to the appropriate 
obstruction calculus; I note the case of~$KO$, with $B\mathrm{Spin}$ in the fiber 
of~\eqref{firstext}. With the obstructions in place, the output is instead a \emph{curved} 
coefficient system, defining \emph{gerby} versions of $\cC_{3,4}$.  

The second, explicit construction realizes the $\cC_{3,4}(G;E)$ by symmetry-breaking to the maximal 
torus followed by Weyl descent. It is tempting to posit an elegant
$A$-model interpretation, 
reflecting the boundary condition for $3$D gauge theory built from the flag variety of $G$ and a polar half 
of $E$; this would be a worthy generalization of the main theorem in \cite{gmp}.

\subsection*{Alternative proposals.} My construction is quite different from the proposal of $\caR_3$ 
in (the updated version\footnote{The original version contained some conceptual errors, which were 
repaired in the sequel.} of) \cite {betal}. While that output ought to agree with (the complexification 
of) mine, the difference in methods makes the correspondence obscure. In particular, the reliance 
on the strange duality from  \cite{missing} does not have a clear analogue in my topological approach. 
Meanwhile, some distinctions are easily explained:
\begin{enumerate}\itemsep0ex
\item My definition of $\caR_3$, as the cohomology of a spectrum, applies with integer coefficients. The 
explicit formulas hold after inverting the order $|W|$ of the Weyl group. In $\caR_4$ (for which~\cite{betal} 
do not propose a construction) the formula holds after inverting only the torsion in~$\pi_1G$.  

\item I propose a stronger obstruction for~$\caR_3$: $w_4(E)$ must have a square root $\bar{r}\in H^2(BG;\bZ/2)$ 
\emph{admitting an integral lift $r$}. This represents the difference between the geometric~$\bC^\times$-action 
on $\cC_3$ and the cohomology grading. For connected~$G$, the condition in~\cite{betal} is equivalent, 
minus the integral lift. Without that, the homology $\bZ$-grading 
is obstructed by the (big) 
Bockstein~$B(\bar{r})$. The smallest example is $\mathrm{SU}(2)\times_{\{\pm1\}} \mathrm{SO}(6)$ 
with its standard representation $\bH\otimes\bR^6$ (Theorem~\ref {whenodd}.ii in the Appendix): $\caR_3$ 
may not be $\bZ$-graded without breaking the Koszul sign rule.\footnote{This is remedied in \emph{loc.~cit.}  
by introducing a separate $\bZ/2$-grading, but that only restates the problem.}

\item The $\bZ$-grading can be collapsed to $\bZ/2$ in complex $K$-theory, and my obstruction
can sometimes be loosened for $\caR_4$. Thus, the example in (ii) above is unobstructed for 
Bott-periodic $K$-theory. It is still obstructed for $KO$-theory, or for connective $kU$-theory, 
and so are the Adams operations on $KU$. 

\item I show (Appendix~A) that, for connected~$G$, a root $\bar{r}\in H^2(BG;\bZ/2)$ of $w_4$ admits 
a lift to $\bZ/4$. In this case, a secondary obstruction~$\sigma$ 
to the existence of~$\caR_4$ appears: 
see Theorem~\ref{kobst}. A typical example is $\mathrm{Sp} (odd)\times_{\{\pm1\}} \mathrm{SO}(4k)$ 
with the tensor product of standard representations.
\item This latter, secondary obstruction seems to supplement Witten's originally proposed obstruction 
for $4$-dimensional~$N=2$ gauge theory with coefficients in~$E$, since $\caR_4$ arises from the 
latter theory by evaluating on $S^1\times S^2$.  The paper \cite{betal} seemed focused 
on validating Witten's primary obstruction~\cite{wit}, which indeed is all that $\caR_3$ can detect.
\end{enumerate}

\subsection*{Improvements and future directions.} 
\begin{enumerate}\itemsep0ex
\item \emph{Alternative descriptions.} The results of Theorem~\ref{one} below are likely not optimal. For example, the natural maps
\[
\cC_{3,4}(G;E)\times_{\cC_{3,4}(G;0)} \cC_{3,4}(G;F) \to \cC_{3,4}(G;E\oplus F),
\]
described on ($K$-)homology coefficients in Part(v), are probably isomorphisms, if we interpret the left side as the affinized fibered quotient under 
the Toda group scheme. 
\item \emph{Disconnected groups.} The explicit classification of obstructions pertains to connected groups: 
the calculations are more involved in general. More significantly, \emph{twisted sectors} appear for 
disconnected groups: the identity sector pertains to point defects in the $3D$ gauge theory, while the twisted sectors represent point defects embedded in topological 
't Hooft loops. 

\item \emph{Categorification of the Coulomb branch.} My $K$-theory coefficient systems $\cK_E$ can be replaced by the respective matrix factorization category (encountering the same obstructions). 
D\'evissage equates the $K$-theory of the resulting category with the one here. It would be interesting to compare 
this categorification with the recent one in \cite{harold}.

\item One interesting question already mentioned: can we interpret the Weyl descent construction 
of the chiral rings in terms of boundary conditions for gauge theory with matter in the spirit of \cite{telc}, generalizing  
the theorem in \cite{gmp}? The spaces $\cC_{3,4}$ come with a Lagrangian multi-section, which could be the quantum cohomology 
of the total space of a vector bundle over the flag variety of $G$ with fiber a Lagrangian half of $E$. 

\item Since the first version of the paper was arXived, interest in \emph{Elliptic Coulomb branches} and
their relation to~$5D$ gauge theory has increased: see the recent discussion in~\cite{bw}. 
The method of this paper indicate a definition of Coulomb branches for arbitrary cohomology 
theories, subject to obstruction-removal. For a version of ellitpic cohomology that is orientable 
for the String group (the $\mathrm{Spin}$ group with $p_1/2$ killed), the obstruction is the composite map
\[
BG \xrightarrow{\ E\ } B\mathrm{Sp} \xrightarrow{\ \eta\ } \mathrm{Sp} \to \mathrm{Sp}_{\le 7} 
\] 
the Postnikov truncation of the group to $\pi_7$. The $K$-theoretic obstruction having 
been settled in this paper, there is an additional step to work out, because the 
new obstructing homotopy group $\pi_7$ is connected by $k$-invariants to the lower ones.  
\end{enumerate}

\subsection*{Organization.} The key results of the paper, Theorems~1--3, are stated in Section~1; 
the reader may need to refer to later sections for some details. In \S2 we discuss two topological 
facts that underlie Theorems~1 and 3. Section~3 quickly reviews the construction \cite{bfn} of 
the polarized case. Section~\ref{real} describes the changes required for the general case: most important is 
the twisted reality structure. In Section~5, we re-interpret the topological obstructions and their cancellation 
in terms of the Weyl group and maximal torus. The global description of chiral rings for the 
normalizer $N(H)$ of the maximal torus in~$G$ is described in \S\ref{Lag}, followed by the adjustment 
needed to pass to the non-abelian group~$G$. The Abelianization theorem, applying to Coulomb branches 
whose representations contain the roots, is briefly discussed in \S8. Appendix~A classifies the 
obstructions for all connected groups, and Appendix~B collects some examples and 
 counterexamples to overly optimistic statements.

\subsection*{Acknowledgements.} This long overdue paper was inspired by a conversation with Sam Raskin 
(Spring 2019), who generously shared his ideas on the subject, as explained in correspondence with  
A.~Braverman~\cite{ras}. Discussion of the primary~$2$-torsion obstruction pointed the way 
out of the ``no-go'' example ($\mathrm{SU}_2$ with its standard representation); the source of 
the obstruction, and its secondary version, from the $\eta$-class in Bott periodicity was then 
easy to guess. The original material is now found in the paper \cite{betal}, whose first version 
appeared while the present one was in preparation.  

I would also like to thank A.~Braverman, T.~Dimofte, M.~Finkelberg, H.~Nakajima and 
H.~Williams for enlightening conversations and Kiran Luecke and especially Ben Webster for spotting 
some mistakes and misstatements. The recent paper~\cite{bw} also contains some enlightening developments. 
The results were obtained during the author's association with MSRI in Spring 2019, and early versions 
were presented at conferences in Moscow, Edinburgh and VBAC~2019 and at the Perimeter Institute in 2020, 
with gratitude to the organizers. Completion of this work 
was supported by the Simons Collaboration 
on Global Categorical Symmetries and by a visit at the Aspen Center for Physics, partially funded 
by NSF grant PHY-1607611. The final version of the paper  benefitted from error-correction 
by Anthropic's Claude, who spotted several typos 
in formulas and proofs, as well as supplying helpful expository pointers.


\section{Background and new results}
\subsection{Pure gauge theory.}\label{sphere top}
The ($K$-)homology rings $\caR_{3,4}(G;0)$ of the based loop group $\Omega G$, equivariant with respect to 
the conjugation $G$-action, were analyzed in \cite{bfm}, and related to the Toda integrable systems. 
The double coset stack $G\backslash LG/G$ of the free loop group $LG$, homotopy equivalent to the 
free mapping stack from $S^2$ to $BG$, is a more revealing model for the natural $E_3$ structure on 
the Pontrjagin multiplication. This $E_3$ structure is shown in 
\emph{loc.~cit.}~to define an 
algebraic symplectic form, while the (Hopf) algebra structures over the ground rings $H^*_G, K^*_G$ 
of a point make the underlying spaces $\cC_{3,4}$ into relative abelian groups, which in addition admit 
integrable system structures. 
They are (fibre-wise group completions of) the classical Toda  system and its finite-difference 
version, respectively. These spaces, which we now denote $\cC_{3,4}(G;0)$,  were later 
interpreted in terms of gauge theory in the guise of classifying spaces for categories with 
topological $G$-action, with  Gromov-Witten theory as motivating example: see \cite{telicm} 
or \cite[\S2]{telc} for a summary. 

\subsection{Polarizable matter.}
With a different motivation, a construction of the spaces $\cC_{3,4}(G;E)$ 
for polarized symplectic representations $E=V\oplus V^\vee$ was  provided by Braverman, 
Finkelberg and Nakajima \cite{bfn}, based on earlier ideas of Nakajima \cite{nak}. 
This involves the ($K$-)homology of a \emph{linear space} $L_V$ over $\Omega G$, an 
algebraic fibration in vector spaces. 

The spaces $\cC_{3,4}$ can also be interpreted in $2$-dimensional gauge theory. A polarization 
allows to couple a \emph{mass parameter} to the matter fields, scaling the two 
polar summands by inverse actions of $\bC^\times$: the equivariant parameter becomes the 
\emph{complex mass} of the physics literature. This also defines two topological boundary conditions 
of the $3$D gauge theory: the Gromov-Witten theories of the spaces $V, V^\vee$, regarded as real $G$-Hamiltonian 
manifolds. Each of them defines a regular Lagrangian section 
of $\cC_{3,4}(G;E) $ over the Toda 
base, and their ratio is a rational section of $\cC_{3,4}(G;0)$. This section can be identified with 
the exponentiated derivative $\exp(d\Psi)$ of the GLSM superpotential $\Psi$. While the full 
structure of these physically inspired constructions has not been completely settled, it is proved in  
\cite{telc} that we recover $\cC_{3,4}(G;E)$ from two copies of the Toda space glued after  
a relative shift by this section.

Both constructions require a polarization: without it, we seem to miss the space $L_V$, and 
no gauge-invariant topological boundary conditions of geometric origin for the 
$3$D theory are apparent that would reconstruct $\cC_{3,4}$. 

\subsection{New results.}
The present paper overcomes these obstacles. First, I adapt the construction of \cite{bfn} by exploiting 
a (twisted) reality structure on $L_V$. The topological obstructions pertain to this twist. 
Second, I polarize $E$ after breaking the symmetry to the maximal 
torus, and descend back under the Weyl group. 
The $K$-theoretic version $\cC_4$ suggests the prospect of nice integral presentations, but those are 
not so obvious from my methods. (Integrality may be the shadow of a categorification of 
$\cC_4$.)  
A topological subtlety requires removing two obstructions, to be discussed in detail in \S\ref{obstruction} 
below.

The constructions anchor a definition. In special cases, a one-step Abelianization (\S\ref{abel}) 
gives a clean answer. In general, the answer is determined\footnote {After some loss of torsion information, 
in integral $H_*$.}  by reductions to the maximal torus and to the  
Levi subgroups of semi-simple rank one, a simplification made possible by the freedom 
of the chiral rings over their Toda bases. The uniform algebraic description condenses this in \S\ref{Lag}. Computability seems to distinguish the construction presented here from other approaches.

\subsection{Statements.} Postponing some details of the \emph{construction}, here are the main 
theorems of the paper. Many proofs are repetitions of arguments in \cite{bfn} or \cite{telc}.  
Assume the removal of the first obstruction for $\cC_3$, and also 
of the secondary one for $\cC_4$. We also assume that $G$ is \emph{connected}.\footnote{Mainly, 
to avoid discussing the twisted sectors.}

\begin{maintheorem}\label{one}
There exist $G$-equivariant, $E_2$-multiplicative coefficient systems $\caH_E, \cK_E$ for 
($K$-)homology over the based loop group $\Omega G$, constructible with respect to Bruhat 
stratification, such that:  
\begin{enumerate}\itemsep0ex
\item The equivariant homology $H_*^G(\Omega G;\caH_E)$ is a 
$\bZ$-graded, commutative algebra over $H^*_G$. \\
With rational coefficients, it is a free $H^*_G$-module. 

\item When $\pi_1G$ has no torsion, this applies integrally to $K_*^G(\Omega G;\cK_E)$ 
over $K^0_G$.\\ (See  Remark~\ref{comp}.iii below for the torsion case.) 

\item These rings carry $E_3$ structures defined by Poisson structures  (of homology degree $2$). These are the leading 
terms of non-commutative deformations, constructed by incorporating the loop rotation action. 

\item 
The Toda group schemes $\cC_{3,4}(G;0)$ act $E_3$-compatibly on $\cC_{3,4}(G;E)$.

\item More generally, there are multiplications 
$\caH_E\otimes \caH_F\to \caH_{E\oplus F}, \cK_E\otimes \cK_F\to \cK_{E\oplus F}$, 
compatible with all the listed structure, once the obstructions have been compatibly cancelled. 
\end{enumerate}
\end{maintheorem}

\begin{remark}[Complements to Theorem~\ref{one}]{\ }\label{comp}
\begin{enumerate}\itemsep0ex
\item The coefficient systems are built from $G$-equivariant constructible sheaves of spectra, 
de-suspensions of the cells of $\Omega G$ by stratified virtual linear spaces $R_E\to\Omega G$ (\S\ref
{coulombspectrum}, \S\ref{real}). These linear spaces are only locally defined, and only up 
to suspension by real, oriented (respectively $\mathrm{Spin}$ or $\mathrm{Spin}^c$) vector bundles. However, the ($K$-)homology 
sheaves $\caH_E,\cK_E$ are unambiguous and multiplicative, once the obstructions have been removed. 
\item With the obstructions in place, we can interpret the $\caH_E,\cK_E$ as \emph{curved} constructible 
coefficient systems:
classes in $H^2(X;\bZ/2)$ define curved local systems for ordinary integral homology. 
\item Given $G=\tilde{G}/\pi$, with $\pi$ finite, and a $G$-space $X$, the ring $K_{\tilde{G}}(X)$ 
is graded by the characters of $\pi$. The character pairing makes $\pi$ act by automorphisms of the ring. 
Taking $\pi$ to be the torsion subgroup of $\pi_1G$, the statements of the theorem 
apply to the orbifolded $K$-theory 
$\pi\ltimes K_*^{\tilde G}(\Omega G;\cK_E\otimes\bC)$ over the orbifold Toda base $\pi\ltimes K^0_{\tilde G}\otimes\bC$. 
\end{enumerate}
\end{remark}

The second construction of the chiral rings $\caR_{3,4}(G;E)$ proceeds by breaking the symmetry to the 
maximal torus $H$,  polarizing $E$ as $E_+\oplus E_-$, and then descending back under the Weyl group. We construct the associated chiral rings $\caR(H;E)$ by a tweak of
the method of \cite{telc}. (Specifically, we convert the Lagrangian shift from \emph{loc.~cit.} into a \emph{charge 
conjugation} automorphism, with the same effect of identifying the desired subring.) 
We then modify the Weyl action on $\cC_{3,4}(H;0)$ by suitable Euler class factors (see \S\ref{modweyl}) 
to compensate for the broken symmetry. 
The Weyl quotients of the subrings surviving 
charge conjugation now give the (identity sectors of the) chiral rings for the normalizer $N(H)$. As in the original construction \cite{bfm}, we obtain the chiral ring for $G$ by an additional modification over the (affine) root hyperplanes; see  
\S\ref{execution} for details.

\begin{maintheorem}\label{two}
The spaces $\cC_{3,4}(G;E)$ are built from the affine quotients of the $\cC_{3,4}(H;E)$  under the shifted Weyl 
action by explicit correction (\S\ref{execution}) on the (affine) root hyperplanes.
\end{maintheorem}

\begin{remark}
It is worth explaining the added complexity in Theorem~\ref{two}, versus the polarized case. 
The Euler class of the index bundle of $E$ defines a rational\footnote
{In the sense of algebraic geometry, not rational homotopy...} torsor over the Toda spaces~$\cC_{3,4}(G;0)$. 
Topologically, the torsor stems from the connecting map $\eta: KSp\to \Sigma^3 KO$, which obstructs the 
polarization lifting of~$E$ to~$KU$: cf.~Wood's sequence \eqref{woodsequence} below. This torsor 
has order~$2$, because $2\eta =0$. (More intuitively, the torsor is multiplicative with respect 
to the representation, and $E^{\oplus 2} \cong E\oplus E^\vee$ is polarized, leading 
to a trivialization of the square.)  There is then some gymnastics involved in identifying the 
modification along the root hyperplanes, needed to retrieve the Coulomb branches for $G$ from $H$.  
\end{remark}

The final formula is a reduction to the Cartan subgroup $H$ and Weyl group $W$ for certain representations 
$E$. This is an \emph{exact} formula, \emph{not} the usual, information-losing localization theorem. 

\begin{maintheorem}[Abelianization]\label{three}
After inverting $|W|$ in the ground ring, 
\[
\cC_{3,4}(G; E) \cong \cC_{3,4}(H;E\ominus\frg_\bH)/W,
\] 
as soon as the roots of $\frg$ appear among the weights of $E$.
\end{maintheorem}

\begin{remark}[Complements to Theorem~3]{\ }
\begin{enumerate}\itemsep0ex
\item
Because $\frg$ is real, $E$ will contain $\frg_\bH$ as $H$-subspace, if the roots appear among the weights. 
\item The assumptions in Theorem~\ref{three} are not optimal.
The $\cC_{3,4}$ are determined by a 
collection of multiplicities associated to the $E$-weight hyperplanes \cite[\S5]{telc}. A root 
hyperplane can be abelianized if and only if 
$E\ge \frg_\bH$ in respect of that multiplicity. Thus, 
$\caR_3\left(\mathrm{SU}(2);E\right)$ 
is abelianizable, except for $E=0$, $E=\bH$ (which is obstructed) and $E=\bH^{\oplus 2}$ (which is polarized), 
determining all the chiral rings. 
\item For $\caR_4$, the multiplicity condition in (ii) also applies to the \emph{affine} root hyperplanes. 
Half-integer multiples of the roots (which occur for symplectic groups)  will not abelianize the space 
over the `odd' affine hyperplanes. Thus, the Chiral ring $\caR_4\left(\mathrm{SU}(2);E\right)$ 
does \emph{not} abelianize over $(-1)\in\mathrm{SU}(2)$ for the irreducible representation $E\cong \bC^4$.
\item Over the base, the Chern character may be used to identify $\cC_4$ locally (analytically) with 
spaces $\cC_3$, for appropriate Levi subgroups of $G$. Generic reduction to the root hyperplanes 
determines every complexified $\cC_{3,4}$.
\end{enumerate}
\end{remark}

\section{Two topological facts}
We review two topology constructions. The first one relies on  \emph{Wood's theorem}, 
implicit in Bott's periodicity of real $K$-theory; it underlies the construction of 
the $\caR(G;E)$ by extracting a `quantum square root' in place of the missing classical one. 
This is analogous to the Spinorial square root of the exterior power of a real vector space, 
which is only obstructed by $w_2$, rather than requiring a complex structure. 
The second fact, the stable splitting of a stratified de-suspension of a manifold, underlies 
the Abelianization  of \S\ref{abel}.

\subsection{First result: Wood's theorem.} \label{wood}
Complexifying a real vector bundle defines a morphism of ring spectra $KO\to KU$. 
Wood's theorem \cite{wood} identifies the resulting fibration sequence as the $KO$-module 
extension of $\Sigma^2KO$ by $KO$ classified by $\eta\otimes: \Sigma^2 KO \to \Sigma^1 KO$: 
\begin{equation}\label{woodsequence}
KO \xrightarrow{\bC\otimes_\bR} KU \xrightarrow{\Omega^2(\bH\otimes_\bC)} \Sigma^2 KO.
\end{equation} 
(This is the only interesting extension, since $\mathrm{Ext}^1_{KO}(\Sigma^2 KO,KO) = 
\pi_1 KO = \{0,\eta\}$.) The map $\Omega^2(\bH\otimes_\bC)$ is the double-looping of the 
quaternionization  $ KU\to KSp = \Sigma^4 KO$, $V\mapsto \bH\otimes_\bC V$. It is also the $\Sigma^2$ of the 
forgetful map $KU\to KO$; both times, we have implicitly used Bott periodicity on $KU$. 
A lift of $\bH\otimes_\bC$ is the datum of a (stable) polarization on 
a quaternionic bundle $E$.

\begin{remark}
The result is not difficult: applying~$KR$ to the two-point set $\{\pm\}$, swapped by 
the Real involution, leads to the fibration sequence
\[
KO \longrightarrow KR(\{\pm\}) \longrightarrow {}^\sigma KO,
\]
in which $KR(\{\pm\}) = KU$, while the twisting $\sigma$ of $KO$-theory is the 
suspension of $KR$ by the (reduced) sign representation of the Real involution. A 
Clifford algebra calculation of the relevant Thom isomorphism identifies the twisted
${}^\sigma KO$ with $\Sigma^2KO$. This argument has the  merit of applying 
equivariantly as well. 

This easy argument is not entirely honest, as it implicitly uses key properties of $K$-theory, 
including Bott periodicity, part of which identifies the complex Lagrangian Grassmannian 
$Sp/U \simeq B(U/O)$ with $\Omega Sp$; and this identification 
already describes Wood's sequence.  
\end{remark}

\subsection{Application.} 
A complex representation $E$ of $G$ gives rise to a complex, virtual \emph{index bundle} 
$\Ind_E$ over $\Omega G$, equivariant with respect to $G$-conjugation. 
The fiber of $\Ind_E$ over a free loop $\gamma:S^1\to G$ is the Dirac index of the $E$-bundle 
$E(\gamma) \to \bP^1$ defined by the equatorial transition function $\gamma$. This 
last description is equivariant for the left$\times$right actions of $G\times G$ on the 
free loop group $LG$, and is compatible with loop rotation.\footnote
{Or rather, the double rotation, with its lift to spinors.}  
A stricter, algebraic construction of $\Ind_E$ arises by interpreting the $G$-equivariant 
homotopy type $\Omega G$ as that of the moduli stack of algebraic $G_\bC$-bundles on $\bP^1$. 
The direct image of the associated $E$-bundles along $\bP^1$, after a half-canonical twist, 
leads to a $2$-term complex of coherent sheaves representing $\Ind_E$. 

An important feature of the index bundle,  stemming from its doubly delooped origin, 
is its two-fold additivity; namely, $\Ind_E: \Omega G \to KU$ is a $G$-equivariant 
$E_2$ map. Commutativity progresses by one step, to $E_3$, when we pass to fixed points, specifically 
the (geometric) 
fixed points of the stable homotopy or $K$-theory linearizations, giving $E_3$-compatible 
maps from $(\Sigma^\infty\Omega G_+)^G$ or  $(K\wedge\Omega G)^G$ to $KU^G$. The $E_3$ structure 
on the source is the \emph{sphere topology product} alluded to in \S\ref{sphere top}; this 
incorporates a wrong-way map, which is why stabilization (in the homotopical sense) is needed.  
We can build an analogous map in ordinary homology if we also change the codomain to the  
$G$-equivariant Eilenberg-MacLane spectrum, for instance, by following the index map with the equivariant 
Chern character. 

\begin{remark}
Despite its topological clarity, this construction faces the difficulty that the multiplication 
is not defined strictly: this makes its application to $\cC(G;E)$ problematic. In effect, one must refine 
it to a stratified map. This problem was solved by the construction in~\cite{bfn}, 
which instead interprets $\Omega G$ as the moduli of bundles over the formal disk with doubled origin. 
This variation does not change the topology. 
\end{remark}

A quaternionic structure on $E$ refines the index to a doubly-suspended real structure:
\begin{equation}\label{koreduction}
\Ind_E:\Omega G \to \Omega Sp=\Sigma^2KO. 
\end{equation}
A polarization $E= V\oplus V^\vee$ supplies a lifting of this to $KU$ in Wood's sequence 
\eqref{woodsequence}:
\begin{equation}\label{polarlift}
\Ind_V:\Omega G \to KU, \quad\text{with}\quad\Ind_E = \Omega^2(\bH\otimes_\bC)\circ \Ind_V. 
\end{equation}
This is used in \cite{bfn} to construct the Coulomb branches $\cC(G;E)$ (see \S\ref{pol} for a 
quick refresher). Applying 
the construction to $\Ind_E$ instead of $\Ind_V$ gives the Coulomb branch 
for the symplectically doubled representation $\bH\otimes_\bC E$. 
The following proposition is used 
in extracting the square root of this construction in the absence of a polarization, as we shall 
do in \S\ref{real}.

\begin{proposition}\label{half}
The map $\Ind_E: \Omega G \to \Sigma^2 KO$ can be lifted locally, $G$-equivariantly, to $KU$. 
\end{proposition}

\begin{remark}\label{liftings}{\ }
\begin{enumerate}\itemsep0ex
\item Liftings form a torsor over $KO_G(\Omega G)$; absent a polarization, there is no preferred lift. 
\item Liftings need \emph{not} be additive, let alone $E_2$. In fact, twice-delooping a global 
equivariant $E_2$ lifting would give a stable $G$-polarization of $E$. Complete reducibility 
of representations would lead us to an actual $G$-polarization. 
\end{enumerate}
\end{remark}

We will exploit the homotopy-equivalent \emph{Laurent polynomial subgroup} 
$\Omega^a G\subset\Omega G$. This is the quotient ind-variety $G_\bC\pz/G\bz$; it is 
stratified by $G\bz$-orbits, which are even-dimensional complex vector bundles over 
the $G$-orbits of the one-parameter subgroups in $G$. The latter are the generalized 
flag varieties $G/L$ of $G$, for various Levi subgroups $L$. 

\begin{proof}[Proof of Proposition~\ref{half}]
The obstruction  $\eta\otimes \Ind_E$ to a local lifting lives in $KO^1_L =0$.
\end{proof}

For later use, we record the following. 
\begin{lemma} \label{KOsplitting}
The stratification of $\Omega G$ splits $KO_*^G(\Omega G)$ into a sum of copies of $KO, KSp$ and $KU$. 
\end{lemma}
\begin{proof} The stratification assembles $KO_*^G(\Omega G)$ from copies of the equivariant 
coefficient rings $KO_L$, suspended by even-dimensional complex representations of the various $L$. 
Each of these is a sum of copies of $KO, KSp$ and $KU$. Since $\mathrm{Hom}_{KO}(M,\Sigma N)=0$ 
for all listed $KO$-modules $M, N$, there are no possible $KO$-linear connecting maps in the Gysin 
sequences for the strata and no $KO$-linear extensions.
\end{proof}

\subsection{Second result: a stable splitting.} 
\label{splitsec}
Let $M$ be a manifold equipped with a Morse-Bott function $f$ whose Morse stratification satisfies 
the Whitney conditions.  For Morse functions, it is proved in \cite{nic} that the Whitney conditions 
are ensured by the Smale transversality conditions\footnote {Plus a technical clustering condition 
on the Hessian eigenvalues at critical points \cite[Remark~4.3.4.b]{nic}; this can be met by adjusting the metric at the critical 
points, and carries no topological content.} \cite{smale}. This is believed to hold in the Morse-Bott case also, 
but I do not know a reference. In our application to $\Omega G$ below, the Whitney conditions will 
be ensured by equivariance. 

Recall that Whitney's Condition (A) asserts that the union $N(f)$ of normals to the strata in the tangent 
bundle $TM$ is closed. We form the Thom spectrum of $N(f)$ and desuspend it by the tangent bundle, 
to obtain a spectrum $\Sigma^fM:= \Sigma^{N(f)-TM}M$, sitting between the Spanier-Whitehead 
dual $\Sigma^{-TM}M$ and the suspension spectrum $\Sigma^\infty M$. 
The Morse stratification of $M$ gives a filtration of $\Sigma^fM$, with associated graded 
spectrum a sum of copies of the sphere $\bS$, one for each critical point. In the Morse-Bott case, 
the slices are suspension spectra of the critical manifolds. 

\begin{proposition}\label{splitprop}
For M a Morse function $f$, $\Sigma^fM$ is naturally a sum of copies of $\bS$. In the Morse-Bott case, this 
is the sum of suspension spectra of the critical strata.
\end{proposition}
\begin{remark}
Interesting self-extensions of $\bS$ are precluded by connectivity, explaining the splitting 
in the case of Morse functions; but this argument does not apply in the Morse-Bott case, and in any 
case does not provide a preferred splitting.
\end{remark}
\begin{proof}
We note the geometric splitting of the attaching maps: as we approach a lower stratum from a higher one, 
we find additional vertical directions in $N(f)$ corresponding to the directions of 
approach, which we can use to shoot out the attaching map towards the base-point at $\infty$.  
\end{proof}

The conclusion of Proposition~\ref{splitprop} also applies equivariantly with respect to compact 
group actions \cite{field}. We will use this in the case of the loop group~$\Omega G$, where the Whitney 
property of the strata of~$\Omega^aG$ follows from homogeneity of the stratification 
under the subgroup $G[\![z]\!]$.

\subsection{Abelianization of certain Coulomb branches.} \label{abelannounce}
The $G\bz$-orbits in the subgroup $\Omega^a G\subset \Omega G$ are the descending Morse-Bott 
strata for the $G$-invariant energy functional $f:\Omega G\to \bR$. Proposition~\ref
{splitprop} splits the spectrum $\Sigma^f\Omega^a G$, $G$-equivariantly, into a sum 
of Spanier-Whitehead duals $\Sigma^{-T(G/L)}(G/L)$ of flag varieties~$G/L$; the sum is 
labeled by Weyl orbits of co-characters of~$G$. 

We now apply equivariant $K$-homology and exploit the isomorphism of coefficient rings 
\[
K^0_L = \left(K^0_H\right)^{W_L}, 
\]
for the Cartan subgroup $H\subset G$ and Weyl group $W_L$ of $L$. This converts the $K_G$ group 
to the~$W_G$-invariant part of the sum, over all co-characters, of copies of~$K_H^0$. The ring structure 
may be tracked by the ordinary localization theorem (see \S4), and doing so recovers the computation 
in \cite{bfn} for the adjoint Coulomb branch ($W=W_G$): 
\[
\cC_4(G; \frg_\bH) \cong \cC_4(H;0)/W.
\]
We generalize this in \S\ref{abel}, to an Abelianization theorem
\begin{equation}\label{abelianformula}
\cC_4(G;E) \cong \cC_4(H;E\ominus\frg_\bH)/W,
\end{equation}
under the assumption that $E$ should contain $\frg_\bH$ as an $H$-representation.

\section{Review of the polarized case: the Coulomb spectrum} 
\label{pol}
We recall the construction in \cite{bfn}  of the 
spaces $\cC_{3,4}(G;E)$ for polarized representations $E = V\oplus V^\vee$, before reframing it 
in terms of the equivariant \emph{Coulomb spectrum}~$\Sigma(G;V)$. More details may be found the 
original paper, and a summary in \cite[\S3, \S6]{telc}, from which the paragraphs below are excerpted. 
I will deviate from the sources by incorporating a Spin structure on the disk from the outset; 
while this clutters the notation with factors of $(dz)^{1/2}$, it avoids subsequent redefinitions. 
Assume that our group $G$ is connected; $\pi_0G$   leads to additional orbifolding.

\subsection{The Chiral rings $\caR_{3,4}(G;0)$ \cite{bfm}.} \label{basiccoul} 

The space  $\cC_3(G;0):= \mathrm{Spec}\,H_*^G(\Omega G;\bC)$ is an affine symplectic resolution 
of singularities of the Weyl quotient $T^* H^\vee_\bC/W$. The homology grading represents 
the $\bC^\times$-scaling of the cotangent fibers. When $G$ is simply connected, 
$\mathrm{Spec}\,K^G_*(\Omega G;\bC)$ is also a symplectic manifold, giving an affine 
resolution of $(H_\bC\times H^\vee_\bC)/W$; in general, it has quotient singularities 
under the torsion subgroup $\pi\subset\pi_1G$. To avoid those, we set $G=\tilde{G}/\pi$, 
$H=\tilde{H}/\pi$, and define $\cC_4(G;0)$ as the smooth symplectic orbifold 
$\pi\ltimes\mathrm{Spec}\, K^{\tilde{G}}_*(\Omega G;\bC)$. 

The Hopf algebra structures of $H_*^G(\Omega G), K^G_*(\Omega G)$ over the ground 
rings $H_G^*, K^{\,0}_{G}$ of a point lead to relative abelian group structures 
\begin{equation}\label{intsyst}
\cC_3(G;\mathbf{0})\xrightarrow{\ \ }\frh_\bC/W, \qquad 
	\cC_4(G;\mathbf{0})\xrightarrow{\ \ }\pi\ltimes (\tilde{H}_\bC/W),
\end{equation}
which also define integrable systems: the first is a fiberwise group completion of the classical
Toda system \cite{bf}, the second is its finite-difference version. These groups 
act on all other Chiral rings.

The $G$-equivariant loop multiplication on $\Omega G$ has 
an algebraic counterpart for $\Omega^aG$, by means of the double coset stack 
$G\bz\backslash  G\pz /G\bz$ and the $G\bz\times G\bz$-equivariant correspondence diagram 
\begin{equation}\label{algmult}
\Omega^a G \times G\bz\backslash G\pz 
	\leftarrow G\pz \times_{G\bz} G\pz \rightarrow G\pz. 
\end{equation}
In both cases, the underlying Poisson structure is the leading term of a non-commutative deformation 
over $\bC[h]=H^*(BR)$, respectively over $\bC[q^\pm]=K^0_R$, obtained by incorporating equivariance 
under the loop-rotation ($z$-scaling) circle~$R$. The analogue  
applies to all Coulomb branches below.

\subsection{The polarized case, $E=V\oplus V^\vee$.} \label{LV}
The  $\cC_{3,4}(G;E)$ are the Specs of the $G$-equivariant ($K$-)homologies of a 
\emph{linear space} $L_V\to\Omega^{a} G$: a $G\bz$-equivariant stratified space 
with vector space fibers. 
Namely, the fiber of $L_V$ over a Laurent 
loop $\gamma\in G\pz$ is the kernel of the difference map
\begin{equation}\label{kay}
\left.L_V\right|_\gamma =\mathrm{Ker}\left\{ V[\![z]\!] \oplus V[\![z]\!] 
\xrightarrow{\mathrm{Id}-\gamma} V(\!(z)\!)\right\}\otimes(dz)^{1/2}.
\end{equation}
This complex is equivariant under the left and right actions of $G[\![z]\!]$ on the 
Laurent loop group, simultaneously acting on the respective factors $V[\![z]\!]$, 
and with the left copy alone acting on $V(\!(z)\!)$. Over any finite set of $G\bz$-orbits 
in $\Omega^a G$, projection to either summand $V[\![z]\!](dz)^{1/2}$  embeds $L_V$ therein, with 
bounded co-dimension. Moreover, $L_V$ also contains two sub-bundles of finite co-dimension, 
from a left and a right $z^nV[\![z]\!]$, $n\gg 0$. 

This stratified finiteness allows \cite{bfn} to define the $G$-equivariant Borel-Moore 
($K$-)homologies of $L_V$, after renormalising the homology grading as if $\dim V[\![z]\!]$ 
were zero. The normalised grading is compatible with the multiplication defined by 
the following correspondence diagram on $L_V$, living over the multiplication of 
two loops $\gamma, \delta\in G\pz$ in the correspondence \eqref{algmult}:
\begin{equation}\label{correspondence}
\left.L_V\right|_\gamma \oplus  \left.L_V\right|_\delta \quad\leftarrowtail\quad
	\left.L_V\right|_\gamma {\underset{V[\![z]\!]}{\oplus}}  \left.L_V\right|_\delta
 	\quad\rightarrowtail\quad\left.L_V\right|_{\gamma\cdot\delta};
\end{equation}
the sum in the middle is fibered over the right component of $\left.L_V\right|_\gamma$ 
and the left one of $\left.L_V\right|_\delta$, while the right embedding projects to the outer 
$V[\![z]\!]$ summands. The wrong-way map in homology along the first inclusion is defined 
after quotienting by a common sub-bundle $z^nV[\![z]\!]$, and the result is 
independent of the sufficiently large $n$. 

\begin{remark} \label{halftwist}
The twist in \eqref{kay} by $(dz)^{1/2}$ is relevant to the loop rotation action and the 
non-commutative chiral rings; here, we only need it to make contact with the Dirac index bundle.
\end{remark}

\subsection{The spectrum $\Sigma(G;V)$.} \label{coulombspectrum}
We re-interpret the construction of $\cC_{3,4}$, removing 
infinite dimensions and the consequent 
renormalization of homology degree, by ``subtracting'' the fiber $V\bz(dz)^{1/2}$ over $1\in \Omega^aG$ 
from the linear space $L_V$. That fiber being the largest of all, the cost of the transaction is 
the passage to stable homotopy.

\begin{definition}
The \emph{Coulomb spectrum} $\Sigma(G;V)$ is the de-suspension of the Thom spectrum of  
$L_V$ by the left bundle $V\bz(dz)^{1/2}$. 
\end{definition}
\noindent
This is a $G\bz$-equivariant stratified de-suspension of $\Omega^aG_+$. It generalizes 
the spectrum $\Sigma^f\Omega^aG$ of \S\ref{abelannounce}, which we obtain for the adjoint 
representation $V=\frg$ (except for the half-integral correction $(dz)^{1/2}$ 
to loop rotation). The correspondence diagram \eqref{correspondence} 
defines an $E_2$ multiplication on $\Sigma(G;V)$, compatible with its inclusion in 
the suspension spectrum $\Sigma^\infty\Omega^aG_+$. The latter is the group ring of 
$\Omega G$ over the  sphere $\bS$, 
and we can think of $\Sigma(G;V)$ as a group ring with coefficients. The function rings 
of $\cC_{3,4}(G;E)$ are the $G$-equivariant ($K$-)homologies of 
$\Sigma(G;V)$.  

\subsection{Left versus right.}\label{leftright}
Another version $\Sigma(G;V)_r$ of the Coulomb spectrum is obtained from the right factor 
of $V\bz$. The ``left minus right'' difference of bundles $V\bz(dz)^{1/2}$ over $\Omega^aG$ 
is the index bundle $\Ind_V$ \cite[\S6]{telc}, so that the two versions are related by 
\[
\Sigma(G;V)_r = \Sigma^{\Ind_V}\Sigma(G;V).
\]
The $E_2$ property of the index bundle makes $\Sigma^{\Ind_V}\Omega^a G_+$ into an $E_2$-ring spectrum: 
namely, the twisted group ring $\Omega G\ltimes_{\Ind_V}\bS$, where $\Omega G$ acts on the ring spectrum 
$\bS$ by composing $\Ind_V$ with the delooping $BJ$ of the $J$-homomorphism: 
\begin{equation}\label{Jind}
\Omega G \xrightarrow {\ \Ind_V\ } BU \xrightarrow{BJ} 
B\mathrm{GL}_{1}(\bS).
\end{equation}
If $c_1V\neq 0$, $\Ind_V$ has non-zero rank and $\Omega G$ maps instead to $\mathrm{Pic}(\bS)$, via 
$\bZ\times BU$.   
From \eqref{Jind}, the left and right Coulomb spectra differ by a central extension of $\Omega G$ by $\mathrm{GL}_{1}(\bS)$. 
Factorization through $BU$ in \eqref{Jind} makes the extension  invisible when applied to a complex-oriented 
homology theory. However, when $c_1V\neq0$, the homology grading on the two versions will differ. 

\subsection{Crossed product formulation.} 
Continuing this idea, call $N_V: = L_V\ominus V\bz(dz)^{1/2}$ our  (virtual) normalization of $L_V$, 
and re-interpret $\Sigma(G;V)$ as the crossed product $\Omega G\ltimes_{N_V} \bS$: $N_V$ 
acts via~$BJ$. The jumps of $N_V$ across strata now give a constructible system of coefficients,  
instead of a bundle. 

We use this  to summarize the construction that follows. The spectrum 
$\Sigma(G;E) = \Omega G\ltimes_{N_E} \bS$ leads to the rings 
$\caR_{3,4}(G;\bH\otimes_\bC E)$ for the double of $E$. We will cut this in half 
by observing that the composition 
\[
\Omega G\xrightarrow{\ \mathrm{Ind}_E\ } BU \xrightarrow{\ BJ\ } \mathrm{Pic}(\bS) \to \mathrm{Pic}(KU),
\] 
trivialized by the complex orientation, has a dimensional half, which almost 
factors through the quotient $\Omega Sp = BU/BO$. Namely, the Atiyah-Bott-Shapiro orientation  
trivialize the real $BJ:BO\to \mathrm{Pic}(KU)$ on the cover $BSpin^c\to BO$. The obstruction 
to a doubly-delooped half of $BJ:N_E\to \mathrm{Pic}(KU)$ is the ($G$-equivariant) composition 
\[
BG \xrightarrow{\ E\ } BSp 
	\xrightarrow{\ \eta\ } \Sigma^3 (\bZ\times BO) 
	\to \Sigma^3\cW_B
\]
with an \emph{obstruction spectrum} $\cW_B$ (Definition~\ref{obspectrum} below), 
incorporating the orientation and $\mathrm{Spin}^c$ obstructions $w_1$ and $W_3$. With the 
obstruction lifted, the symplectic structure of $E$ can be used to define $\Omega G\ltimes_{N_E} KU$ 
(for $\caR_4$,  and its homology version for $\caR_3$).


\section{General case: twisted real structure on the linear space}
\label{real}
Subject to obstructions and ambiguities to be discussed, we now construct the coefficient systems 
$\caH_E, \cK_E$ over $\Omega^a G$ substituting for the non-existent stable homotopy incarnation $\Sigma(G;V)$.

\subsection{Stratified polarization of $\,\Ind_E$.}
\label{stratpol}
We can think of the linear space $N_E:=L_E\ominus E\bz(dz)^{1/2}$ as a constructible lift of the index
$\Ind_E$ to $\bZ\times BU$, with respect to Wood's sequence \eqref{woodsequence} over the strata 
of $\Omega^a G$. Stratum-by stratum, 
\begin{equation}\label{localift}
\Omega^2(\bH\otimes_\bC)\circ N_E = \Ind_E. 
\end{equation}
Indeed, near any one-parameter subgroup $z^\gamma$, $\gamma\in\frh$, $E$ breaks up as $E_+\oplus E_0\oplus E_-$ according to $\gamma$-eigenvalue and $E_-$ polarizes the complementary representation $E\ominus E_0$ of $L$. 
Realizing $\Ind_E$ as the $\bar{\partial}$-cohomology of the bundle $E(\gamma)\otimes \sqrt{K}$ over $\bP^1$, 
$H^0$ comes from $E_+$ and $H^1$ from $E_-$. ($E_0$ does not contribute.) Comparing with 
\[
z^\gamma: E\bz(dz)^{1/2} \to E\pz/E\bz(dz)^{1/2} 
\]
identifies $N_E$ with $(-H^1)$, and Serre duality plus the quaternionic structure on $E$ 
give the anti-linear identification of $H^0$ and $H^1$ required for \eqref{localift}. This refined 
interpretation sees the left side of \eqref{localift} as a degeneration of the right as a 
constructible coefficient system: moving transversally to the $\gamma$-stratum deforms the $\bar{\partial}$-operator on $\bP^1$,   
providing an extension class which converts the symplectic double of $N_E$ into $\Ind_E$. 

\subsection{Real structures.}
A real structure on $N_E$ would provide a de-suspension of $\Omega G$ reaching half-way to $\Sigma(G;E)$, 
analogous to the effect of a polar half $V$. Now, on neighborhoods $U$ of $G\bz$-orbits, $\eta$ can 
be trivialized (Proposition~\ref {half}) and a second, continuous $G$-invariant local lift $S_E$ of 
$\Ind_E$ gives a stable real structure on $N_E\ominus S_E$. The jump in the fibers comes from $N_E$. 
Denote the underlying (locally defined) real linear space by~$R_E$; writing 
\[
\bC\otimes R_E = N_E\ominus S_E
\]
makes the suspension $\Sigma^{R_E}U$ into a real version of the de-suspension of $\Sigma(G;E)$ by $S_E$: 
\[
\Sigma^{\bC R_E}U_+ = \left.\Sigma^{-S_E}\Sigma(G;E)\right|_U. 
\]

\subsection*{Example: Polarized case.}\label{polexample}
When $E=V\oplus V^\vee$, we can take $\Ind_{V^\vee}$ for $S_E$. Denote by underlines the $\sqrt{K}$-twists 
of the associated bundles on $\bP^1$. Then, 
\[\begin{split}
N_E- S_E &= -H^1\left(\underline{E}\right) -H^0\left(\underline{V}^\vee\right) 
	+H^1\left(\underline{V}^\vee\right) \\
	&= -H^1\left(\underline{V}\right) -H^0\left(\underline{V}^\vee\right) 
	= -H^1\left(\underline{V}\right) -H^1\left(\underline{V}\right)^\vee;
\end{split}\]
the underlying real space is $-H^1\left(\underline{V}\right) = L_V\ominus V\bz(dz)^{1/2}$, and 
$\Sigma^{R_E}\Omega^aG_+ = \Sigma(G;V)$.

\subsection{Obstruction theory.}\label{obstruction}
Applying $H^G_*,K_*^G$ to $\Sigma^{R_E}U_+$ should give the constructible coefficient systems 
$\caH_E,\cK_E$ defining our chiral rings $\caR_{3,4}(G;E)$. However, we meet two problems:
\begin{enumerate}\itemsep0ex
\item The local Thom spectra $\Sigma^{R_E}U_+$ depend on the 
auxiliary local section $S_E$, resulting 
in a suspension ambiguity by a general $G$-equivariant $KO$-class. 
\item Even if $S_E$ is globally defined, there is no Pontrjagin product, if $S_E$ is not multiplicative. 
\end{enumerate}
The local ambiguity in (i) becomes a global obstruction if $\eta\otimes \Ind_E\neq 0\in KO^1_G(\Omega G)$. 
The product in (ii) runs into the group extension of $\Omega G$ by $KO$, pulled back 
from Wood's sequence \eqref{woodsequence} by $\Ind_E$; this gives a \emph{projective} action of 
$\Omega G$ on the sphere $\bS$, pre-empting the crossed product ring  $\Omega G\ltimes_{R_E} \bS$ 
of \S\ref{leftright}. An $E_2$ splitting of the $KO$-extension is equivalent to a polarization of $E$ 
(Remark~\ref{liftings}.i). Without it, there is no \emph{stable homotopy} chiral ring. 
Fortunately, the obstructions to building $\caR_{3,4}$ are much  milder. We will encounter  

\begin{itemize}\itemsep0ex
\item the \emph{primary obstruction}: $w_4(E)\in H^4(BG;\bZ/2)$ modulo squares $r^2$, $r\in H^2(BG;\bZ)$; 
\item the \emph{secondary obstruction}: $B\sigma\in H^6(BG;\bZ)$ as below.  
\end{itemize}
After trivializing~$w_4-r^2$ by a cochain $t\in C^3(BG;\bZ/2)$, 
\[
\sigma:=Sq^2 t + \eta\cdot w_4\in H^5(BG;\bZ/2),
\] 
and $B$ is the integral Bockstein. The vanishing of $Sq^2(r^2) = Sq^3Sq^1(r)$, and the trivialization 
$Sq^2(w_4) = \delta(\eta\cdot w_4)$ on~$B{Sp}$, convert 
$Sq^2t+\eta w_4$ into a co-cycle, well-defined up to 
$Sq^2H^3(BG;\bZ/2)$. These two classes can be packaged into the \emph{obstruction spectrum} defined 
below, but we first state and explain the

\begin{theorem}[Primary and Secondary obstructions]\label{summaryprop}{\ }
\begin{enumerate}\itemsep0ex
\item Construction of $\caR_3(G;E)$ requires lifting the 
\emph{primary obstruction}. Thereupon, the choices of $\cC_3$ form a torsor over $H^3(BG;\bZ/2)$. 
\item Constructing $\caR_4$ also requires removing the \emph{secondary obstruction $B\sigma$}; 
thereupon, the choices form a torsor over the $2$-torsion in $H^5(BG;\bZ)$. 
\item The two obstructions combine to a map $BG\to \Sigma^3\cW_B$ (Definition~\ref{obspectrum} below). They both vanish when $G$ is connected without symplectic factors (Theorems~\ref{whenodd}, \ref{kobst}). The ambiguities combine to a map $BG\to \Sigma^2\cW_B$. 
\end{enumerate}
\end{theorem}

\begin{remark}[Complements] {\ }\label{obsremark}
\begin{enumerate}\itemsep0ex
\item With complex (rather than integer) coefficients, the choices in $H^3$ affect $\caR_3$ only via 
their Bockstein images in (the $2$-torsion in) $H^4(BG;\bZ)$, so their effect vanishes for connected $G$. 
\item When $\pi_0G,\pi_1G$ have no $2$-torsion, $H^5(BG;\bZ/2)=0$, so \emph{a fortiori} $\sigma$ and $B\sigma$ 
vanish. 
\item The class $\sigma$ itself obstructs the $KO$-version of $\caR_4$; see \S\ref{koobst} below.
\item If we are willing to collapse the homology grading mod~$2$, the primary obstruction may be 
removed by a square root $\underline{r}$ of $w_4$ in $H^2(BG;\bZ/2)$; see~\S\ref{otherobst}.
\item When $E=V\oplus V^\vee$, the universal identity 
$c_2(E) = 2c_2(V)-c_1^2(V)$ cancels the obstructions. Conceptually, the lift of $E:BG\to BSp$ by  $V:BG\to BU$ 
kills the obstruction source $\eta$. 
\item In both cases, the $E_2$ multiplication improves to $E_3$ on \emph{equivariant} ($K$-)homologies:  
the transgressive nature of the obstructions tracks their 
cancellation on the space of maps $S^2\to BG$.
\end{enumerate}
\end{remark}

\begin{remark}[Ambiguities] \label{choices}
The choices in Theorem~\ref{summaryprop} are meaningful in TQFT. The rings $\caR_{3,4}(G;E)$ are expected 
to be the vector spaces associated to the sphere in a $3$-dimensional gauge theory with group $G$ and 
matter fields in~$E$. The pure gauge theory, without matter fields, can be precisely, 
if incompletely,\footnote{In the sense that not all $3$-dimensional operations are defined.} 
defined as the \emph{sphere topology} of the stack $BG$. To a closed surface $M$, this assigns 
the ($K$-)homology of the moduli $\mathrm{Bun}_G(M)$ of holomorphic $G$-bundles on $M$ (which is 
homotopy equivalent to the mapping space of $M\to BG$). 
These theories admit discrete twists 
by $\Sigma^2\cW_B(BG)$. 

A twist integrates over~$M$ to a class in $\cW_B\left(\mathrm{Bun}_G(M)\right)$, defining a twisting 
for ($K$-)homology over~$\mathrm{Bun}_G(M)$, modifying the groups. Nonetheless, the effect 
is easily concealed for \emph{connected~$G$} and \emph{complex coefficients}. For~$\cC_3$, a 
class $w\in H^3(BG;\bZ/2)$ transgresses trivially over surfaces: this is because $w$ lifts to a 
Bockstein image from $H^2(BG;\bZ/2)$. The change is visible on the \emph{integral} ring structure 
(see \S\ref{so3workout}), but disappears in complex cohomology, because the image of $w$ vanishes 
in $\bC^\times$ coefficients. 
Next, a typical trivialization of~$B\sigma$ comes from one of~$\sigma$ (see Theorem~\ref{kobst}), 
confining the ambiguity to the Bockstein image of $H^4(BG;\bZ/2)$, which vanishes.  
\end{remark}

\begin{definition}\label{obspectrum} 
The \emph{obstruction spectrum} $\cW_B$ is the co-fiber of the map $BSpin^c\to \bZ\times BO$, with its natural deloopings.
\end{definition}
\noindent
This $\cW_B$ is closely related to the bottom $3$ layer truncation of $KO$; its three homotopy groups are 
\[
\pi_0\cW_B =\bZ, \pi_1=\bZ/2, \pi_3=\bZ, 
\]
having shifted $\pi_2(BO)=\bZ/2$ to~$\pi_3=\bZ$ by the big Bockstein~$B$. The $k$-invariant $B\circ Sq^2$ links the top two groups into a spectrum 
$\cW_B'$; the latter is placed in $\cW_B$ over the base $\bZ$ by the lift of $Sq^2$ in $[H\bZ; \Sigma\cW'_B]$, seen 
to be unique from the sequence 
\begin{equation}\label{obskinvariant}
[H\bZ; \Sigma^4\bZ] =0 \xrightarrow{} 
[H\bZ; \Sigma\cW'_B] \xrightarrow{\ \sim\ } \langle Sq^2\rangle = [H\bZ; \Sigma^2\bZ/2] \xrightarrow{BSq^2} [H\bZ; \Sigma^5\bZ] 
\end{equation}
in which $BSq^2\circ Sq^2 =0$.

\begin{proof}[Proof of Theorem~\ref{summaryprop}, semi-simple $G$]
The Thom isomorphism in homology reconciles  
the local choices in \S\ref{obstruction}.i, provided we reduce the 
structure group $BO$ to $BSO$ in the 
sequence \eqref{woodsequence} over $\Omega G$: 
\[
BO\to BU \to \Sigma^2(\bZ\times BO).
\] 
Reduction to $B{Spin}^c$  accomplishes the same for complex $K$-theory. 
Doubly delooping such reductions leads to an $E_2$ multiplication.
Executing this meets the double-deloopings of the orientation obstruction $w_1$ (for $\caR_3$) and 
the additional $Spin^c$ 
obstruction~$W_3 = Bw_2$ (for $\caR_4$), jointly classified by the composition
\begin{equation}\label{ob2}
BG\xrightarrow{\ E\ }  BSp \xrightarrow{\eta\otimes} \Sigma^3BO \to \Sigma^3\cW'_B.
\end{equation}
Cancelling the bottom component $w_4(E)$ as $\delta t$ leaves the remaining obstruction to defining $\caR_4$ as the Bockstein image 
$B\sigma = B(Sq^2(t)+\eta w_4)\in H^6(BG;\bZ)$. 
\end{proof}

\begin{proof}[Proof for reductive $G$.]\label{reduct}
Polarized representations of a reductive group may have odd $c_2$, leading to a nonzero obstruction in 
$\cW'_B$: such is the symplectic double of the standard representation of $\mathrm{U}(1)$. The hidden problem is that the fiber of the middle map $\eta\otimes:BSp\to \Sigma^3BO$ in \eqref{ob2} is 
$BSU$, rather than $BU$, preventing lifts by a polar half with $c_1\neq 0$. 

Abandoning the $\bZ$ at the base of $\cW_B$ carried a cost. 
While $BSp$ has no interesting maps to $\Sigma^3\bZ$, the trivial map there out of $BG$ has self-homotopies  
classified by $H^2(BG;\bZ)$. (This is zero for semi-simple $G$.) The effect of $r\in H^2(BG;\bZ)$ is to shift $w_4$ by $r^2 \pmod{2}$. Its effect on the full 
$\cW_B$-obstruction can then be found from the sequence \eqref{obskinvariant}.  
\end{proof}

\begin{remark}
The freedom to remove the $w_4$ obstruction is also explained by the formula 
\[
c_2(E\oplus L\oplus L^\vee) = c_2(E)-c_1(L)^2
\] 
for a  one-dimensional representation $L$. The polarized summand $L\oplus L^\vee$ is killed by $\eta$ 
and does not change the obstruction $\eta\otimes E$.
\end{remark}

\subsection{Some basic properties of $\caR_{3,4}$.} 
To prepare the algebraic construction of the chiral rings in \S\ref{Lag}, we now check 
their evenness and freedom over the Toda base.   

\begin{proposition}[Parity check] {\ }
\begin{enumerate}\itemsep0ex
\item If $w_4(E)=\bar{r}^2$, with $\bar{r}\in H^2(BG;\bZ/2)$, then $R_E$ may be chosen with even-dimensional 
(real) fibers. 
\item If $\bar{r}$ has an integral lift, these dimensions can be chosen to define an even $\bZ$-grading for the 
Pontrjagin product on $\caR_3(G;E)$.
\end{enumerate}
\end{proposition}

\begin{proof}
With notation as in \S\ref{stratpol},
\[
\left.\dim_\bR R_E\right|_{z^\gamma} = -\sum\nolimits_\nu \langle \nu | \gamma\rangle - \dim_\bC S_E, 
\]
with $\nu$ over $\gamma$-positive weights of $E$ with their multiplicities. Call $q$ the quadratic 
form 
on the Lie algebra associated to $c_2$. As $k=k^2\pmod{2}$, we may switch to the sum of squares 
$\langle \nu |\gamma\rangle^2$, which evaluates to $q(\gamma)$. For two co-weights $\gamma,\gamma_0$ 
in the same component of $\Omega G$, 
\[
q(\gamma) - q(\gamma_0) = \partial^2q(\gamma_0,\gamma-\gamma_0) + 
	q(\gamma-\gamma_0),
\]
with the associated bilinear form $\partial^2 q$. The latter is even, because $c_2$ is a 
sum of integral squares. Under our assumptions, $q$ takes even values on co-roots, since 
$\bar{r}=0$ there. Adjusting $\dim_\bC S_E$ on components can then render $\dim R_E$ even, which proves (i). 

With an integral lift $r$, we may specifically choose $\dim S_E = r(\gamma)$ to settle Part~(ii).   
\end{proof}

\begin{corollary}[Freedom over the Toda base]
The $\caR_3(G;E)\otimes \bQ$ is free as a module over the base $H^G_*(\text{point})\otimes \bQ$. The same applies 
integrally to $\caR_4$, when $\pi_1G$ is torsion-free. (See Remark~\ref{comp}.iii for torsion in $\pi_1$).
\end{corollary}

\begin{proof}[Proof of the Corollary.]
The Bruhat stratification filters those groups with free module subquotients in 
even degrees, ruling out connecting differentials.
\end{proof}

We conclude with two variations on the obstruction theory $\cW_B$. 

\subsection{$KO$-obstruction.}\label{koobst}
As mentioned, $\cW_B$ can be built using the big Bockstein $B$ from the $3$-layer truncation 
$ko_{< 3}$: a spectrum $\cW$ with $\pi_0=\bZ$,  
$\pi_1=\pi_2=\bZ/2$, with top 
$k$-invariant $Sq^2$. This $\cW$ carries the secondary obstruction $\sigma$ to real spinnability.
One  handle on $\sigma$ comes from the formula, shadowing the Adem relation $Sq^2Sq^2 = Sq^3Sq^1$,
\begin{equation}\label{sigmarel}
Sq^2\sigma = Sq^3\frac{c_2-r^2}{2}
\end{equation}
This stems from the restriction of $ko_{< 3}$ 
to $2\bZ\subset\pi_0\, ko$, which is built by stacking $\Sigma^2\bZ/2$ over $2\bZ\times\Sigma\bZ/2$ with 
$k$-invariant $(x,y)\mapsto Sq^3(x/2)+Sq^2y$.

\subsection{Abandoning the homology grading.}\label{otherobst}
The classification in the Appendix shows that, for \emph
{connected}~$G$, the obstruction discussed in  
\cite{betal} is equivalent to the weaker requirement of a mod~$2$ square root of $w_4$. Allowing that, 
only Case (i) of Theorem~\ref {whenodd} is obstructed. Now, the bottom two layers $\bZ,\Sigma\bZ/2$ 
of~$ko$ may be collapsed to $\bZ/2,\Sigma\bZ/2$, and for $\cC_3$, we gain the freedom of 
cancelling $w_4(E)$ by a square from  $H^2(BG;\bZ/2)$. However, the bottom $\bZ$ represents the 
grading in homology: exploiting $H^2(BG;\bZ/2)$ collapses the homology $\bZ$-grading to $\bZ/2$. 

The obstruction to the grading is mild, and is only manifested at the second delooping: the Koszul sign rule is broken in the Pontrjagin multiplication.


\section{Obstructions in terms of the Weyl group}
\label{allweyl}
I now review the topological obstructions in terms of the maximal torus of $G$ and 
its normalizer. Their \emph{transgression} to the free loop group of $BG$ is relevant to a modified Weyl group 
action on the Toda space, to be used in the explicit construction of the chiral rings 
in the next section. One benefit clarifies the elusive secondary obstruction $\sigma$: its transgression 
has a well-defined leading term in ordinary cohomology.

\subsection*{Notation:} 
$H\subset G$ denotes the maximal torus, $H^\vee$ 
the dual torus, ${}^{2}H^\vee\subset H^\vee$ the subgroup of $2$-torsion points, 
$N(H)$ the normalizer of $H$ in $G$, $W$ the Weyl group, $W_{\mathrm{aff}}$ the affine Weyl group $W\ltimes \Lambda^\vee$, 
with the weight lattice $\Lambda$. 
Roots and coroots will be denoted by $\alpha$ and $h_\alpha$, general weights of representations by $\nu$. 
Typical elements of $H$ and $\frh$ will be denoted by $x$ and $\xi$, respectively. 
We will also use a generic regular 
element $\xi_0\in \frh$, splitting $E= E_+\oplus E_0 \oplus E_-$ into positive, negative and zero 
weight spaces.\footnote{$E_0$ must come entirely from polarizable summands in~$E$; see~Appendix~A.} 
Finally, for~$w\in W$ and a map $f:\frh \to H^\vee$, $w[f]$ will denote the map $\xi\mapsto w\left(f(w^{-1}\xi)\right)$.

\subsection{Transgressed obstructions.}\label{weylobst}
The evaluation map $S^1\times LBG \to BG$ on the loop space defines a transgression,\footnote{Using the 
bounding Spin structure on the circle} with the adjoint action of $G$ on itself on the right-hand side: 
\[
\tau: KSp^0_G \xrightarrow{} KSp^{-1}_G(G)
\]
Further restriction takes us to $KSp^{-1}_{N(H)}(H)$. Our two obstructions from $\eta\cdot E$ lead to
\[
\tau(w_4) \in H^3_{N(H)}(H;\bZ/2),\qquad
\tau(\sigma) \in H^4_{N(H)}(H;\bZ/2);
\]
the second one is contingent on a cancellation of the first, but part of 
it turns out to be independently defined (Remark~\ref {welldefined} below). We pare them 
down to their leading terms in the Leray spectral sequence 
$H^p_W\left(H^q(BH\times H)\right) \Rightarrow H^{p+q}_{N(H)}(H)$, to obtain classes in Weyl group cohomology: 

\begin{theorem}{\ }\label{weylcocycles}
\begin{enumerate}\itemsep0ex
\item The Leray leading term of $\tau(w_4)$ is a class in 
$H^2_W\left(\Lambda/2\right) \subset H^3_W(H;\bZ/2)$, 
represented by
\[
(u,v)\mapsto c(u,v) := \sum\nolimits_{\substack{\nu>0\\ v\nu<0\\ uv\nu>0}} uv\nu \in \Lambda/2.
\]

\item The Leray leading term of $\tau(\sigma)$ is a class in $H^1_W\left(\Lambda^{\otimes 2}/2\right)\subset H^4_W(BH\times H;\bZ/2)$. It is uniquely defined, modulo the $Sq^2H^3$ ambiguity in $\sigma$, 
and is represented by 
\[
w\mapsto s^{\otimes 2}(w):= \sum_{\substack{\nu>0 \\ w\nu<0}} w\nu\otimes w\nu \in \Lambda^{\otimes 2}/2. 
\]
\item The Leray leading term of the big Bockstein $B\tau(\sigma)\in H^2_W\left(\Lambda^{\otimes 2}\right)$ is represented by
\[
(u,v)\mapsto d(u,v):= \sum\nolimits_
	{\substack{\nu<0\\ v\nu>0\\ uv\nu<0}} uv\nu\otimes 
uv\nu \in \Lambda^{\otimes 2}.
\]
\item A trivialization of $w_4$, respectively $(w_4,\sigma)$ or $(w_4, B\sigma)$ leads to one for the Weyl 
cocycles $c$, respectively $s^{\otimes 2}$ and $d$. 
\end{enumerate}
\end{theorem}
\noindent
The cohomology degrees disambiguate the two copies, $H^1(H)$ and $H^2(BH)$, of $\Lambda$ in the formulas.

\begin{remark} \label{welldefined}
{\ } \begin{enumerate}\itemsep0ex
\item The obstruction $\sigma$ is only defined upon trivializing $w_4$. However, the cancellations 
we find below imply the vanishing of the $H^{0,1}_W$ components of 
$\tau(w_4)$, and define $\tau(\sigma)$ 
unambiguously over the $3$-skeleton of $BW$ even when $w_4\neq0$. 
\item An example with $w_4=0$ but $s^{\otimes 2}\neq0$ is given 
in \S\ref{s2not0}. In particular, $\sigma\neq 0$ there.
\item The Weyl cocycles have more trivializations than the topological obstructions for~$G$.  
Cancellations coming from~$G$ \emph{must} be used in building $\cC_{3,4}$. A sample erroneous  cancellation, preventing the construction of $\cC_3$, 
is given in  \S\ref{wrongweyl}.
\item When the obstructions vanish, the cocycles $c$ and $d$ can be cancelled by Weyl $1$-cochains of 
the form $\sum\nu$ and $\sum \nu\otimes \nu$. This follows from 
judicious choices of orientations (Remark~\ref{nar}).
\item Formula~(i) also applies to $\tau(w_4+\bar{r}^2)$, with $\bar{r}\in H^2(BG;\bZ/2)$: indeed, as 
$\deg\tau(\bar{r}) = 1$,
\[
\tau(\bar{r}^2) = \tau\left(Sq^2\bar{r}\right) = 
Sq^2\tau(\bar{r}) =0.
\] 
\end{enumerate}
\end{remark}

\subsection{Preparation: symplectic case.} 
The formulae are functorial under restriction from the representation $G\to \mathrm{Sp}(E)$; 
this reduces the verification to $\mathrm{Sp}(m)$ with its standard representation $\bH^m$. 
Specifically, $H$ maps to the standard Cartan subgroup 
$H_{\mathrm{Sp}}\cong \mathrm{U}(1)^{\times m}\subset \mathrm{Sp}(E)$, giving a compatible map of normalizers, Weyl groups and 
dual map $H_{\mathrm{Sp}}^\vee \to H^\vee$: 
\[
W_{\mathrm{Sp}}=S_m\wr \bZ/2, \quad N\left(H_{\mathrm{Sp}}\right) = 
	S_m\wr \mathrm{Pin}^-(2), \quad 
	\{\nu| \nu>0\} = \{\mathbf{e}_j\},
	\quad H_{\mathrm{Sp}}\cong H_{\mathrm{Sp}} ^\vee = \bC^m/\langle \mathbf{e}_j\rangle.
\] 
Each $\mathbf{e}_j$ flips signs under the $j\textsuperscript{th}$ factor of $\bZ/2$ in $W$, and the last isomorphism identifies co-roots with fundamental weights.  

Identify $\Lambda/2$ with the $2$-torsion in $H_{\mathrm{Sp}}$: then,  $c\in H^2_W(\Lambda/2)$ defines the normalizer extension 
\[
1\to \mathrm{Spin}(2)^m \to S_m\wr \mathrm{Pin}^-(2) \to S_m\wr \bZ/2 \to 1,
\]
or more precisely, its Tits reduction to a $2$-torsion extension \cite{tits}. 
An explicit account of Tits extensions can be found in \cite[\S5]{dw}, but it is easy to 
check the symplectic case directly: the permutation group, and the $\bZ/4$-extensions of the $\bZ/2$ 
factors, embed jointly in $\mathrm{Sp}$, and we can compare this extension with the formula for $c$.

\begin{proof}[Proof of Theorem \ref{weylcocycles}]
For Part (i), we show a stronger statement: over the $4$-skeleton of $BW$,
\begin{equation}\label{chern2}
w_4(E) = \sum_{\nu>0} \nu^2 + [c] \in H^4(BN;\bZ/2) 
	\mod{H^4_W(\bZ/2)},
\end{equation}
with $[c]\in H^2_W(\Lambda/2)$ as in Part~(i), but reinterpreted as a class in $H^2_W\left(H^2(BH;\bZ/2)\right)$. 
The $\nu$-sum  is unambiguous below $H^4_W(\bZ/2)$: the first Leray differential is $d_3$, so 
the~$\nu$ are defined over the $2$-skeleton of $BW$  with ambiguity in $H^2_W(\bZ/2)$. 
Their squares are then well-defined mod~$2$ and $H^4_W(\bZ/2)$. 

Before checking \eqref{chern2} for $\mathrm{Sp}$, let us see how it implies Part~(i). First, 
$\tau Sq^2\nu = Sq^2(\tau\nu)=0$; furthermore, $\tau\equiv 0$ after its restriction 
to the subspace $BN(H)\subset LBN(H)$ of constant loops. The $H^0_W$ and $H^1_W$ components 
of~$\tau(w_4)$ thus vanish, and our leading term for $\tau(w_4)$ is indeed in~$H^2_W\left(\Lambda/2\right)$, 
and arises from $[c]$ under the identification $\tau:H^2(BH) \xrightarrow{\sim} H^1(H)$. 

Returning to~$\mathrm{Sp}(m)$, the wreath product decomposition of its normalizer reduces our claim 
to the group $\mathrm{Pin}^-(2) \subset \mathrm{SU}(2)$, with Weyl group $\bmu_2$. Call $x$ the 
generator of $H^1_{\bmu_2}(\bZ/2)$ and $\nu$ that of $\Lambda = H^2(B\mathrm{Spin}_2;\bZ)$. 
The extension class in $\mathrm{Pin}^-(2)$ leads to a non-zero Leray differential $d_3(\nu) = x^3 \in H^3_{\bmu_2}(\bZ/2)$. 
Then, $d_3(\nu^2)=0$, so that~$\nu^2$ defines a co-cycle over 
the $4$-skeleton of~$BW$ whose extension to a co-chain on $B \mathrm{Pin}^-(2)$ satisfies
\[
\delta\nu^2 = \delta Sq^2\nu = Sq^2\delta\nu = Sq^2 x^3 = x^5.
\]  
Now, $\nu^2$ agrees with $w_4 = c_2\bmod{2}$ on $B\mathrm{Spin}(2)$, so that $\nu^2-w_4$ defines 
a class in $H^2_{\bmu_2}(\Lambda/2)$ on the $E_2$ page with Leray differential 
\[
d_3(\nu^2-w_4) = x^5 = d_3(x^2\nu).
\]  
This means that  $\nu^2-w_4 = x^2\nu = c \in H^2_{\bmu_2}(\Lambda/2)$, confirming \eqref{chern2}. 

For Part (ii), note that $\tau(w_4)$ starts in $H^{>1}_W$, so that $\tau(\eta\cdot E)$ has well-defined 
cohomological leading terms in $H^{0,1}_W$ of $H^\bullet(BH\times H;\bZ/2)$. These represent the 
respective components of $\tau(\sigma)$, whenever restricting to a $(G, E)$ where $\sigma$ is defined: 
ambiguities from $H^3(BG;\bZ/2)$ transgress with leading term in 
$H^1_W\left(H^1(H)\right)$, 
whose $Sq^2$ vanishes on the $1$-skeleton of $BW$, and will not alter these leading terms. Moreover, 
as we are transgressing (in $KSp$) over the bounding spin circle,\footnote{We continue to the $2$-sphere, for 
$\Omega G$.} the integral of $\eta$ vanishes and then so does the $H^0_W$ component of $\tau(\sigma)$.

Returning to $\mathrm{Sp}(m)$ for the calculation, note that the cocycle $s^{\otimes 2}$ is strictly 
invariant under $S_m \subset W$, and the same applies to $\tau(\sigma)$; so we are again reduced 
to a check for $\mathrm{SU}(2)$, when $s^{\otimes 2}$ is the generator of $H^1_{\bmu_2}(\bZ/2)$. 
With $x,\nu,\tau(\nu)$ as before, we show that over the $3$-skeleton of $B\bmu_2$,
\[
\tau(\sigma) = x\nu\tau(\nu)\in 
	H^1_{\bmu_2}\left(H^3(BH\times H;\bZ/2)\right).
\] 
In light of the earlier justification, we ignore the (nonexistent) trivializing co-chain for $w_4$ for the calculation, 
and set $\sigma=\eta\cdot w_4$; then, over the $4$-skeleton of $B\bmu_2$, 
\[
\begin{split}
\delta(\eta\cdot w_4) &= Sq^2(w_4),\quad\text{ so that} \\  
\delta(\tau(\sigma)) &= Sq^2(\tau(w_4)) 
	= Sq^2\left(x^2\tau(\nu) + O(x^3)\right) 
= x^4\tau(\nu) + O(x^5).
\end{split}
\]
The only class in the Leray sequence cancelling this (via $d_3$) is the advertised class $x\nu\cdot\tau(\nu)$.   

Part (iii) is a direct calculation of the coboundary on group cochains. 

Part (iv) follows by tracking a topological trivialization under the transgression $\tau$. 
\end{proof}

\subsection{Modified Weyl action.} \label{modweyl}
We now interpret the transgressed obstructions in terms of  modified \emph{rational} actions of the Weyl group 
on the Toda spaces for $H$, lifting the natural actions 
on the bases $\frh$ and $H$. Let $E_+\subset E$ be the sum of positive weight spaces in the 
standard representation of $\mathrm{Sp}(m)$. To an element $w\in W$, assign the rational sections   
\begin{equation}\label{cocycles}
\begin{split}
\chi_w: \frh\to H^\vee,\quad & \xi\mapsto 
\prod_ {\substack{\nu>0 \\ w\nu<0}} 
\langle w\nu|\xi\rangle^{w\nu}, \\
\kappa_w: H \to H^\vee,
\quad &x\mapsto  \prod_ {\substack{\nu>0 \\ w\nu<0}} 
\left(1-x^{-w\nu}\right)^{w\nu},
\end{split}
\end{equation}
and modify the action of $w$ on Toda spaces by $\chi$- and $\kappa$ shifts, as follows:
\[
\begin{split}
(\xi,h) &\mapsto \left(w\xi,\,\chi_w(w\xi)\cdot wh\right), \\
	(x,h) &\mapsto \left(wx,\,\kappa_w(wx)\cdot wh\right).
\end{split}
\]
This Weyl action will not quite close. To see this, identify $\Lambda/2\cong {}^2H^\vee$ and 
$\mathrm{Hom}\left(H;H^\vee\right)\cong \Lambda^{\otimes 2}$, and recall the Weyl cocycles $c$ and $d$ 
from Theorem~\ref {weylcocycles}. Each product $cd(u,v)$ can then be viewed as an affine map $H\to H^\vee$.

\begin{proposition}[Projective obstructions] 
\begin{equation}\label{extensions}
\begin{split}
\delta\chi = c: \quad \chi_{uv} &= \chi_u\cdot u\left[\chi_v\right] \cdot c(u,v), \\
\delta\kappa = c\cdot d: \quad \kappa_{uv} &= \kappa_u\cdot u\left[\kappa_v\right] \cdot c(u,v)\cdot d(u,v)\end{split}
\end{equation}
\end{proposition}

\begin{proof}
We check the formula for $\chi$. Denoting by $\varphi_\nu = 
\langle\nu|\xi\rangle^\nu$, we have
\[
	\varphi_{-\nu} = (-1)^\nu \varphi_\nu^{-1}, \quad
	\varphi_{w\nu} = w\left[\varphi_\nu\right] \text{ for }w\in W, 
\]
and then 
\begin{align*}
 \chi_u\cdot u\left[\chi_v\right]  
	&= \prod_{\substack{\nu>0\\ u\nu<0}} u\varphi_\nu \cdot 
	\prod_{\substack{\nu>0\\ v\nu<0}} uv\varphi_\nu = \prod_{\substack{v\nu>0 \\ 
	uv\nu<0}} uv\varphi_\nu \cdot \prod_{\substack{\nu>0 \\ v\nu<0 \\ uv\nu>0}} uv\varphi_\nu
\cdot\prod_{\substack{\nu>0\\ v\nu<0 \\ uv\nu<0}} uv\varphi_\nu \\
	&= \prod_{\substack{v\nu>0 \\ uv\nu<0}} uv\varphi_\nu\cdot
	\prod_{\substack{\nu<0 \\ v\nu>0 \\ uv\nu<0}} uv\varphi_\nu^{-1}
	\cdot\prod_{\substack{\nu<0\\ v\nu>0\\ uv\nu<0}} (-1)^{uv\nu}
\cdot
	\prod_{\substack{\nu>0\\ v\nu<0 \\ uv\nu<0}} uv\varphi_\nu \\
	& = \prod_
	{\substack{\nu>0 \\ v\nu>0 \\ uv\nu<0}} uv\varphi_\nu \cdot \prod_{\substack{\nu>0 \\v\nu<0 \\ 
	uv\nu<0}} uv\varphi_\nu
	\cdot\prod_{\substack{\nu<0\\ v\nu>0\\ uv\nu<0}} (-1)^{uv\nu} = 
	\prod_{\substack{\nu>0\\ uv\nu<0}} uv\varphi_\nu = \chi_{uv}\cdot c(u,v)^{-1};
\end{align*} 
in the second line, we changed the sign of the label $\nu$ in the middle factor.
The check for $\kappa$ uses the relation
$\psi_{-\nu} = (-x^\nu)^{-\nu} \cdot \psi_\nu^{-1}$ for the corresponding factors, leading to the extra factor $d(u,v)$.
\end{proof}

\begin{corollary}
A cancellation of the obstructions $w_4$ and $\sigma$ leads to a cancellation of the projective 
Weyl cocycles~\eqref{extensions}, by means of $1$-cochains in ${}^2H^\vee$ for $c$, and by the Bockstein 
of an element in $\Lambda^{\otimes 2}/2$ for $d$. \qed
\end{corollary}

\subsection{Interpretation.} In Weyl group cohomology, the cocycle $\chi$ represents the connecting image 
(generalized Bockstein) of the $H$-equivariant Euler class of $\mathrm{Ind}(E)$ for the fibration sequence 
\[
BSO \rightarrowtail \bZ\times BU \twoheadrightarrow \bZ\times\Omega Sp^{w_4}.
\]
($Sp^{w_4}\to Sp$ is the covering killing $w_4$.) Specifically, $\chi$ is the Euler class of the co-boundary 
of the ($H$-equivariant) lift $\Ind_{E_+}$ of $\Ind_E$ to $BU$. 
The analogue holds for $\kappa$ and the $K$-theory Euler classes, with the extension
\[
BSpin^c \rightarrowtail \bZ\times BU 
	\twoheadrightarrow \bZ\times \Omega Sp^{w_4,B(\eta w_4)}
\]
(after we killed,  in $BSp$,  both $w_4$ and the Bockstein of the newly created class $\eta\cdot w_4$).   

\subsection{Obstruction-removal.}\label{nar}
To reconcile the Weyl and topological obstructions and their removal, consider, for 
$w\in W$, the ($\Lambda$-valued) $H$-equivariant Dirac index $D(w)$ over $\bP^1$ defined by 
\[
\langle D(w)|\gamma\rangle := \Ind_{wE_+(\gamma)}\ominus 
\Ind_{E_+(\gamma)}.
\] 
The $D(w)$ have a real structure, as they stem from the difference of two polarizations of $E$ (cf.~\S\ref{real}). 
They carry $H$-actions with no invariant lines, and are thus always individually ($K$-) orientable. The Weyl cocycles $c$ and $Bs^{\otimes 2}$ 
stem from inconsistent orientation choices. 

More precisely, $(w,\gamma)\mapsto \langle D(w)|\gamma\rangle$ defines a class in 
$\bH^2_{W_{\mathrm{aff}}}(kO_H)$, with leading component in $\bH^1_W(\Lambda\otimes kO_H)$ because of the linearity in $\gamma$. Therein,  

\begin{itemize}\itemsep0ex
\item The leading term in $H^1_W(\Lambda)$ is the Bockstein of the Weyl-invariant element $w^H_2(E_+)\in \Lambda/2$. This vanishes when $w_4(E)=0$; otherwise,
it leads to the appearance of some $2$-torsion points in the calculus of the next section, e.g. Remark~\ref{boundaryth}.
\item The orientation obstruction $w_1\in \bZ/2$ takes us to a class in $H^2_W(\Lambda/2)$, 
represented by the Weyl co-cycle $c$; vanishing of $w_4(E)$ allows a choice of consistent 
orientations. 
\item The Spin obstruction, $w_2^H\in \Lambda/2$, leads to the class in $H^1_W(\Lambda^{\otimes 2}/2)$  
represented by $s^{\otimes 2}$; when $B\sigma$ vanishes, we can choose consistent $\mathrm{Spin}^c$ orientations. 
\item The Weyl cocycles $\chi$ and $\kappa$ are the ($K$-theory) Euler classes of the $D(w)$.
\end{itemize}

\begin{remark}[Spin case]
Even when $c$ and $s^{\otimes 2}$ are trivializable, the individual $D(w)$ need \emph{not} be spinnable. 
Their $KO$-Euler classes (Spin determinants)  are then naturally sections of $H^\vee$-bundles of 
order~$2$ over $H$, rather than maps to $H^\vee$ (Remark~\ref{kol} below). Obstructions are cancelled by coherent identifications between the Weyl-transformed 
bundles and their determinant sections.
\end{remark}

\section{Weyl descent for Coulomb branches: the normalizer}\label{Lag}
We move to the explicit algebraic construction of the $\cC_{3,4}(G;E)$ extending the polarized case of 
\cite[Theorem~2]{telc}. One tweak is the use of a \emph{charge conjugation symmetry C}, 
replacing the vertical shift of \emph{loc.~cit.}; this effects a good clean-up of signs in the formulae. 
The main result, Theorem~2, is restated in precise form in the next section. Due to the technicalities 
and formulas, it will help to include a brief navigational 
chart:
\begin{itemize}\itemsep0ex
\item \S\ref{polarcase} recalls the GLSM construction for the polarized case; 
\item \S\ref{modif} reformulates that in terms of $C$; 
\item Remark~\ref{boundaryth} (optional and partly speculative) discusses the relation to boundary theories;
\item \S\ref{genoutline} outlines the general construction;
\item \S\ref{normH} describes the (identity\footnote{Coulomb branches for 
disconnected groups have twisted sectors, 
from the components of the group.} component of) the chiral ring for the normalizer  $N(H)$;
\item Finally, \S\ref{execution} completes the construction of 
$\caR(G;E)$.
\end{itemize}

\subsection{Review of the polarized case.}
\label{polarcase} 
We ``couple a complex mass term'' to $E= V\oplus V^\vee$: we scale $V$ and $V^\vee$ under 
opposite actions of $S^1$, and adjoin the respective equivariant 
parameters, $\mu\in H^2(BS^1)$ or 
$e^{\pm\mu} = m^{\pm 1} \in K^0_{S^1}$, to the bases of the respective Toda systems. We then introduce 
the rational, Lagrangian \emph{Euler sections} 
$\vep_V, \lambda_V$ of the Toda projections,  
\begin{equation}\label{polareuler}
\begin{split}
\vep_V &:\xi\in \frh \mapsto \prod\nolimits_{\nu} 
	\left(\langle\nu|\xi\rangle+\mu\right)^\nu, \\
	\lambda_V &: x\in H\mapsto\prod\nolimits_{\nu} 
\left(1-m^{-1}x^{-\nu}\right)^\nu,
\end{split}
\end{equation}
with $\nu$ ranging over the weights of $V$ (multiplicities included). In parsing 
these formulae, note the two uses of $\nu$: as an infinitesimal character of $H$, and as a lattice vector 
in $\frh^\vee = \mathrm{Lie}(H^\vee)$ (so that $a^\nu$ is a point of $H^\vee$ when $a\in \bC^\times$). Remarks~\ref
{superpotential}--\ref{kol} below review the meaning of these formulae.
\begin{theorem*}[\cite{telc}, Theorems 1 and 2]
The (massive) chiral rings 
$\caR_3(G;E)[\mu]$ and $\caR_4(G;E)\left[m^\pm\right]$ comprise those regular functions on the Toda 
space (with $\mu$, $m^\pm$ adjoined) which remain regular under vertical shift by 
$\vep_V$, respectively $\lambda_V$. Specializing to $\mu=0$, respectively $m=1$, gives the massless versions. 
\end{theorem*}
\noindent

\begin{remark}[Mass terms] 
The auxiliary step of adding $\mu, m$ is not necessary when the group $G$ contains a circle acting on $V$ 
with strictly positive weights. Absent such a circle, we enlarge $G$ to include the 'mass circle' $M$,
scaling $V$ and $V^\vee$ with opposite weights. From the Coulomb branch for $G\times M$, we obtain the massive 
Coulomb branch for $G$ by collapsing the Toda fiber $M^\vee_\bC$.
\end{remark}

\begin{remark}\label{sectioneval}
The affine spaces $\cC(G;E)$ submerge onto their Toda bases away from co-dimension~$2$; 
they are thus determined by their regular local sections. For instance, following  
the description in \cite{bfm}, regular are those sections of $\cC(G;0)$ which pull back to sections $s$ of $\cC(H;0)$  satisfying the evaluation conditions $\exp(h_\alpha)\circ s = 1$ 
over the root hyperplanes $\alpha=0$, or respectively  $e^\alpha =1$. To construct $\cC(G;E)$, we simply include $\vep_V$ or $\lambda_V$, along with their 
Toda group translates, as regular sections, and affinize the resulting space. 
\end{remark}

\subsection{Charge conjugation C.}\label{modif}
It will prove useful to modify the Lagrangian shift, combining it with the automorphism of order $2$
of the Toda space $\cC(G;0)$ which dualizes the $G$-representations 
and simultaneously changes orientation on the $2$-sphere. The former effects inversion on the Toda base, $\xi\leftrightarrow -\xi$, or $x\leftrightarrow x^{-1}$; 
the latter is inversion in $\Omega G$, thus on the Toda fibers. This gives 
\[
C: (\xi,h) \leftrightarrow \left(-\xi,h^{-1}\right), \qquad
	(x,h) \leftrightarrow \left( x^{-1},h^{-1}\right).
\] 
Clearly, functions in $\caR_{3,4}(G;E)$ are also characterized by regularity under the \emph{modified charge conjugations $C_V$} combining the shift with $C$: 

\begin{equation}\label{almostinvol}
C_V: (\xi,h)\mapsto \left(-\xi, \vep_V^{-1}\cdot h^{-1}\right), \qquad
(x,h) \mapsto \left(x^{-1}, \lambda_V^{-1}\cdot h^{-1}\right).
\end{equation}
The $C_V$ square to the vertical shifts by the sections
\begin{equation}\label{invertiblebdry}
\xi\mapsto \prod\nolimits_\nu(-1)^\nu,\qquad\text{respectively}\quad 
x\mapsto \prod\nolimits_\nu (-x^\nu)^{-\nu}, 
\end{equation}
which define automorphisms of $\cC_{3,4}(G;E)$. The first gives a sign automorphism of $\caR_3(G;E)$, 
grading it by the class $w_2(V) \in H^2(BG;\bZ/2)$. The second acts on $\caR_4(G;E)$ by tensoring with the 
(possibly graded) line bundle on $\Omega G$ transgressed from $\frac{1}{2}c_2(E)$. 

\begin{remark}[Interpretation: boundary theories] 
\label{boundaryth}
The sections in \eqref{invertiblebdry} correspond to \emph {invertible topological boundary theories} 
for topological $3$D $G$-gauge theory, with homology and $K$-linear coefficients, respectively. 
Algebraically, these stem from topological actions of $G$ on the categories $\mathrm{Vect}$ of vector spaces and $K$-Mod 
of modules for $K$-theory, respectively. The first action is determined by the 
$\bmu_2$-extension of $G$ classified by $w_2^G(V)$; similarly, 
$c_2(E)/2$ defines an extension 
of $G$ by the group of complex lines. The interpretation of boundary conditions as Lagrangian sections of the Toda space is 
the one outlined in \cite{telicm}.

Alternatively for $\caR_4$, the $2$-extension of $G$ by 
$\bC^\times$ defined by $c_2(E)/2$ 
gives an invertible boundary theory for $4$D topological gauge theory. Transgressing it to a central 
$\bC^\times$-extension of the free loop group $LG$ defines a topological action of the latter 
on $\mathrm{Vect}$,  and thus an invertible boundary 
theory for $3$D $LG$-gauge theory. 
The spaces $\cC_4$ may thus be considered as Coulomb branches for the $LG$; 
and the second section in \eqref{invertiblebdry} is the Lagrangian associated to this boundary theory.

By construction, $\vep_V$, $\lambda_V$ become \emph{regular} Lagrangian sections of the Coulomb branches 
$\cC_{3,4}$. The shifted automorphism $C_V$ of $\cC_3$ cycles between the unit section, $\vep_V^{-1}$ and their translates 
by the section $\prod\nolimits_\nu(-1)^\nu$. One expects these sections to correspond to the $2$-dimensional 
boundary theories for $E/G$ associated to the two polar halves $V, V^\vee$. The preferred 
``unit'' section is a feature (or bug) of the construction of the Coulomb branch, not intrinsic to it.  
When $w_2(V)\neq 0$, there 
is a choice in gauging the linear Sigma model, related to the role of Spin structures in 
defining Floer and symplectic cohomology; this is reflected in the torsion sections $\prod\nolimits_\nu(-1)^\nu$.  

As defined, $C_V$ has infinite order on  $\cC_4$. It can be modified to have finite order, 
but this involves a torsor of order $2$ over the Toda space (Remark~\ref{kol} below). 
The latter becomes necessary for a $KO$ version of $\cC_4$. I shall not work out the 
details for $KO$ here.
\end{remark}

\subsection{The general case in outline.}\label{genoutline}
Polarize $E$, after first reducing the symmetry to $H$. Weyl symmetry must be initially broken.    
The construction by charge conjugation is then quotiented out 
by the modified Weyl group action of \S\ref{allweyl}. The abelian branches $\cC_{3,4}(H;0)$ produce the (identity sectors of the) 
$\cC\left(N(H);E\right)$. An adjustment of multiplicities on the root hyperplanes, matching the 
construction of $H_*^G(\Omega G)$  in \cite{bfm}, yields the chiral rings $\caR_{3,4}(G;E)$.

\subsection{Euler Lagrangians.}
Recall the splitting $E= E_+\oplus E_0 \oplus E_-$; the $H$-invariant part $E_0$ of $E$ will not contribute 
in what follows. Consider the maps to $H^\vee$
\begin{equation}\label{eulerdef}
\vep_+:\xi\in \frh \mapsto \prod_{\nu>0} \langle\nu|\xi\rangle^\nu, \qquad
	\lambda_+ : x\in H\mapsto\prod_{\nu>0} 
\left(1-x^{-\nu}\right)^\nu.
\end{equation}
The graphs of these maps are Lagrangian (see below) and their closures are smooth over the complement 
of their indeterminacy loci on the bases $\frh, H$; these have co-dimension~$2$.

\begin{remark}[Polarized case]\label{polarreduction}
When $E=V\oplus V^\vee$, we can include the mass parameter $\mu$ as an extra line in the Toda base, 
and choose $\xi_0$ along it. Then, $E_+=V, E_-= V^\vee$, recovering the `massive Lagrangians' 
$\vep_V,\lambda_V$ of \eqref{polareuler}. In this case, $W$-invariance of the splitting reduces the modified Weyl 
action of \S5 to the natural one, recovering the construction for the polarized case. 
\end{remark}

\begin{remark}[GLSM]\label{superpotential}
$\vep_+$ and $\lambda_+$ are the exponentiated differentials of the superpotentials for the mirror 
of the GLSM for $E_+$ by $H$:
\[
\xi \mapsto \Psi(\xi)=\mathrm{Tr}_{E_+}\left(\xi(\log \xi -1)\right), \quad
x\mapsto \Psi(x)=\mathrm{Tr}_{E_+}\mathrm{Li}_2(x).
\]
\end{remark}

\begin{remark}[Index interpretation]\label{kl}
A co-weight $\gamma$ of $H$ defines a character $\exp(2\pi\gamma)$ of $H^\vee$, as well as an $E_+$-fiber 
bundle $E_+^\gamma\to \bP^1$. Then, $\exp(2\pi\gamma)\circ\vep_+$ 
and $\exp(2\pi\gamma)\circ\lambda_+$ 
are the equivariant Euler classes, in cohomology and $K$-theory, 
of the Dirac index of $E_+^\gamma$ over $\bP^1$. 
\end{remark}

\begin{remark}[Spin orientation]\label{kol}
More canonical than $\lambda_+$ is the Euler class in $KO$-theory,\footnote{This would be needed for the $KO$-version of $\cC_4$.} 
\[
\lambda^O_+: x\in H\mapsto\prod_{\nu>0} 
\left(x^{\nu/2} - x^{-\nu/2}\right)^\nu.
\]
which is a section of the $H^\vee$-torsor of order $2$ over $H$, defined by the bilinear form 
\[
\sum_{\nu>0} \nu\otimes \nu: (\pi_1H)^{\otimes 2} \to \bZ/2.
\] 
This represents a real central extension of $H$ by $H^\vee$; it can be non-trivial even in unobstructed 
situations, such as for $\mathrm{SU}(6)$ acting on $\Lambda^3\bC^6$. 
\end{remark}

\subsection{Construction of $\cC^1_{3,4}(N(H);E)$.}
\label{normH}

Recall the modified Weyl action of \S5, adjusted to cancel the projective cocycle. We have the

\begin{proposition}[Modified charge conjugation] \label{regweyl}
The modified Weyl action commutes with the following rational ``$C_+$ automorphism'' of Toda spaces: 
\[
\begin{split}
(\xi, h) &\mapsto \left(-\xi, \vep_+^{-1}(\xi)\cdot h^{-1}\right),\\
(x,h) &\mapsto \left(x^{-1}, \lambda_+^{-1}(x)\cdot h^{-1}\right).
\end{split}
\]
\end{proposition}
\begin{proof}
We write out the check for $\chi$ ($\kappa$ is analogous):
\[
\begin{split}
\left(\xi, h\right)& \xmapsto {\ C_+\ } \left(-\xi,\: \vep_+^{-1}(\xi)\cdot h^{-1}\right) 
	\xmapsto{\ Weyl\ }  \left(-w\xi,\: \chi_w(-w\xi) \cdot 
	w\!\left(\vep_+^{-1}(\xi)\right)\cdot wh^{-1}\right) \\
\left(\xi, h\right)& \xmapsto{\ Weyl\ }
	 \left(w\xi,\: \chi_w(w\xi)\cdot wh\right) \xmapsto{\ C_+\ } 
	\left(-w\xi,\: \vep_+^{-1}(w\xi)\cdot \chi_w(w\xi)^{-1} 	wh^{-1}\right),
\end{split}
\]
and the equality of right-hand sides amounts to the (easily checked) relations 
\[
\begin{split}
\chi(\xi)\chi(-\xi) &= \delta\vep_+(\xi), \\ 
\kappa(x)\kappa(x^{-1}) &= \delta\lambda_+(x),
\end{split}
\]
or more precisely, 
\begin{equation}\label{vepchikappa}
\begin{split}
\chi_w(\xi)\chi_w(-\xi) &= w\left[\vep_+\right](\xi)\cdot\vep_+^{-1}(\xi), \\
	\kappa_w(x)\kappa_w(x^{-1}) &= w\left[\lambda_+\right](x)\cdot\lambda^{-1}_+(x).
\end{split}
\end{equation}
The identities persist after correcting $\chi$ and $\kappa$ by $1$-cochains, valued in ${}^2H^\vee$ or 
in $\Lambda^{\otimes 2}$.  
\end{proof}

\begin{theorem}\label{normalcoulomb}
The identity components of $\caR_{3,4}(N(H);E)$ comprise precisely the modified-Weyl invariant functions 
on $\frh\times H^\vee$, respectively $H\times H^\vee$ which remain regular under the automorphisms $C_+$.
\end{theorem}

\begin{proof}
The topological construction proceeds from the linear spaces $H^0\left(E^\gamma(-1)\right)$ over 
$\Omega H/N(H)$: we subtract the Dirac index bundles of local polar halves of $E$ and define real 
structures and ($K$-theory) orientations on the resulting spaces. Interpreting `local' as `on the Weyl 
cover' and using the polarization $E_+$ leads to 
\[
H^0(E^\gamma(-1)) \ominus H^0\left(E^\gamma_+(-1)\right) \oplus H^1\left(E^\gamma_+(-1)\right) =
	H^0\left(E^\gamma_-(-1)\right)\oplus H^0\left(E^\gamma_-(-1)\right)^\vee.
\] 
The rings $\caR_{3,4}(H;E)$ are then spanned, over their Toda bases, by the Fourier modes labeled by 
the co-weights of $H$, coupled to the ($K$-theory) Euler classes of the (underlying real) bundles $H^0\left(E_-^\gamma\right)$: 
\begin{equation}\label{span}
\begin{split}
(\xi, h)&\mapsto h^\gamma\cdot \prod\nolimits_{\substack
	{\langle\nu|\gamma\rangle >0\\ \nu<0}} 
\langle\nu | \xi\rangle^{\langle\nu|\gamma\rangle}, \\
(x,h)&\mapsto h^\gamma\cdot \prod\nolimits_{\substack
	{\langle\nu|\gamma\rangle >0\\ \nu<0}} 
\left(1- x^{-\nu}\right)^{\langle\nu|\gamma\rangle}. 
\end{split}
\end{equation}
We can see that the inclusion of these Euler classes 
\begin{itemize}\itemsep0ex
\item is compatible with the modified Weyl action of \S\ref{modweyl}, and
\item converts the charge conjugation $C$  into the $C_+$ Proposition~\ref{regweyl}.
\end{itemize}
The characterization of the spans of \eqref{span} given in Theorem~\ref{normalcoulomb} should now be clear.
\end{proof}


\section{The non-abelian Coulomb branches}

\subsection{From $\caR^1\left(N(H)\right)$ to $\caR(G)$.}
\label{execution} 
As in the polarized case, $\caR(G;E)$ is a super-ring of $\caR^1\left(N(H)\right)$, contained in its 
localization away form the (affine) walls. The poles over the walls are controlled by restricting 
the local sections of~$\cC\left(N(H);E\right)$ over the Toda base on which the functions in~$\caR(G)$ 
are required to evaluate regularly. Flatness of~$\caR$ over the Toda base reduces us to Levi subgroups 
of semi-simple rank one, which \'etale-cover the generic parts of the walls. 

We `straighten' the Weyl action  near a wall by the sections~$r,q$ of Lemma~\ref{sqroots} below,  
then correct the space to match the starting point $\cC_{3,4}(G;0)$, rather than $\cC_{3,4}(H;0)/W$. 
(The obstructions concern the existence of~$r$ and~$q$.) The recipe is cleaner over the \emph
{interesting} walls where~$\caR$ does \emph{not} abelianize (as in the next section): single-valued 
test sections suffice, while double-valued ones are otherwise needed. I describe the base 
cases~$\mathrm{SU}(2),\mathrm{SO}(3)$ and $\mathrm{U}(2)$ shortly, but first, an explicit description of the outcome may help:

\begin{itemize}\itemsep0ex
\item When no multiples of a root~$\alpha$ appear among the weights of~$E$, the modified Weyl 
action is regular generically on the $\alpha$- or $e^\alpha$-wall, as is $C_+$; and we then  graft 
the original$\cC(G;0)$ onto $\cC(N(H);E)$ there. 
\item When an integral multiple of~$\alpha$ appears as a weight,\footnote{It must then have even multiplicity, 
because the respective irreducible representation is orthogonal.}  the outcome is equivalent to removing 
a factor of~$\langle\alpha | \xi\rangle^{2\alpha}$ from~$\vep_+$ and the~$\chi$ factors, or of~$(1-x^{-\alpha})^{2\alpha}$ 
from their $K$-theory counterparts: this is the abelianization theorem in the next section. 
\item When only half-integer multiples\footnote{This can only happen for the long roots in type $C$.} 
of $\alpha$ occur in $E_+$, but their coefficient sum exceeds~$1$, we correct $\cC_3$ as before; 
but for $\cC_4$, we only remove a factor of $(1-x^{-\alpha/2})^{2\alpha}$, and graft the original 
$\cC_4(G;0)$ on the wall $x^{\alpha/2} = -1$. 
\item There is a final intermediate case, where exactly two copies of~$\alpha/2$ (and no other multiples) 
appear in~$E$. (This happens when $E$ contains the double of the standard representation 
of a symplectic factor of $G$.) 
\end{itemize}

\subsection{The case $\mathrm{Pin}^{-}(2)\subset\mathrm{SU}(2)$ .}
\label{su2workout}
The rings~$\caR_3$ for these groups depend only on an integer~$N\ge 0$: the half-sum of positive 
weights in $E$. (They are obstructed when the sum is odd.) With the same~$E$, passing from~$\mathrm{Pin}^{-}(2)$ 
to~$\mathrm{SU}(2)$ can be described as a shift $N\to (N-2)$ (plus a tweak for $N=0$), the formal effect 
of subtracting a copy of $\frg_\bH$ from $E$. Explicit formulas are in~\eqref{caRPin} and \eqref{caRSU2} 
below, but do mind the \emph{modified coordinate}~$Z$, with respect to the original 
coordinates~$z,\tau$ on $H^\vee, \frh^*$. 

\begin{remark}
I will omit the analogous story for $\caR_4$, which incorporates 
a second integer $N'$, the half-sum of 
\emph{even} weights, controlling the correction over the central point $(-\mathrm{Id})$ (the affine wall). Note 
that any pair of integers $N\ge N' \ge 0$ with even~$N'$ may be obtained from a polarizable representation. 
The topological construction makes it clear that the rings depend only on the ($H$- or $K$-theoretic) Euler 
classes of the index bundles~$E^\gamma$ over the loops $z^\gamma\in \Omega\mathrm{SU}(2)$, and those are determined by the numbers~$N$, respectively $N,N'$ 
of~\S\ref{su2workout}. This confirms the formulas below in general.
\end{remark}
The modified Weyl action is 
\[
z\leftrightarrow (-1)^N\tau^{2N}z^{-1},\qquad \tau\leftrightarrow (-\tau)
\]
and the Weyl invariant part $\caR_3(H;E)^W$ of the chiral ring for $H$ is generated by
\[
\tau^2,\quad u_N:=\tau(z-(-1)^N\tau^{2N}z^{-1}), 
	\quad v_N:= z+(-1)^N\tau^{2N}z^{-1},
\]
giving  
\begin{equation}\label{caRPin}
\caR_3(H;E)^W =\bC[\tau^2, u_N, v_N], \qquad 
u_N^2 = \tau^2v_N^2-4(-1)^N\tau^{2N+2}.
\end{equation}
In the modified coordinate $Z=\tau^{-N}z$ on $H^\vee$, the Weyl action  becomes the natural one: 
$Z\leftrightarrow Z^{-1}$. Then,
\[
u_N = \tau^{N+1}(Z-(-1)^NZ^{-1}), \quad v_N= \tau^N(Z+(-1)^NZ^{-1});
\]
whereas the $\mathrm{SU}(2)$ version of the ring, using $Z$ as standard coordinate, is~\cite{bfm, bfn} 
\begin{equation}\label{caRSU2}
\caR_3(\mathrm{SU}(2);E) = \bC[\tau^2,\: \tau^{-2}u_N, \:
\tau^{-2} (v_N- 2\delta_{0N})].
\end{equation} 

\subsection{Sections for $\mathrm{Pin}^-(2)$.} 
\label{pinsect}
We  pass from~\eqref{caRPin} to~\eqref{caRSU2} by requiring the  functions in $\caR_3
(\mathrm{SU}(2);E)$ be regular only on certain sections $u_N(\tau),v_N(\tau)$ of 
$\cC_3(\mathrm{Pin}^-(2);E)$, with prescribed behavior at~$\tau=0$. They are defined from 
(formal) Laurent paths $\tau\mapsto z(\tau)\in H^\vee$. The  relation in~\eqref{caRPin} shows,  
for the valuation $Val_\tau$ at~$\tau=0$,
\begin{equation}
Val_\tau(u_N) \ge 1 +  \min\left(Val_\tau v_N, N\right), \quad
	Val_\tau(v_N) \ge \min\left(Val_\tau u_N, N+1\right) -1
\end{equation}
with equality if $Val_\tau v_N \neq N$, $Val_\tau(u_N)\neq N+1$, or if we can freely 
scale~$u_N$, or~$v_N$. 

For generic~$v_N(0)$, this means that 
$Val_\tau(v_N) = 0$, $Val_\tau(u_N)=1$, and
\[
Val_\tau\left(P(\tau^2, v_N) + u_N\cdot Q(\tau^2, v_N)\right) = \min\left(Val_\tau(P), 1+ Val_\tau(Q)\right)
\]
with cancellations prevented by the parity mismatch; so the ring~\eqref{caRPin} cannot be 
enlarged within $\caR[\tau^{-2}]$ while preserving regularity on generic~$z(\tau)$. 

\begin{proposition}\label{sings}
The functions in $\caR_3(H;E)^W[\tau^{-2}]$ which are regular on general sections meeting the 
restrictions below are generated over $\caR_3(H;E)^W$ by $u_N/\tau^2$, $(v_N - 2\delta_{0N})/\tau^2$.  
\begin{enumerate}
\itemsep0ex
\item When $N=0$: $v_0(\tau) = 2 \pmod{\tau^2}$; equivalently $z(0) = Z(0) =1$. \\ This gives $Val_\tau(u_0) = 2$. 
\item When $N=1$: $Val_\tau(v_1) =2$; equivalently, $z(\tau) = \pm\tau \pmod{\tau^2}$, or $Z(0) = \pm 1$. 
\\ This gives $u_1(\tau) = \pm2\tau^2 \pmod{\tau^4}$.   
\item When $N\ge 2$: $v_N(\tau) \in O(\tau^2)$; equivalently, $2\le Val_\tau(z) \le 2N-2$, or $|Val_\tau(Z)| \le N-2$. 
\\ The constraints are imposed at the extremes of $Val_\tau$. 
Generically, $Val_\tau(u_N) =3$.
\end{enumerate}
\end{proposition}
\begin{proof}
Valuations make the claim clear in case (iii). Case~(i) relies on exploiting the unrestricted 
values of~$u_0''(0), v_0''(0)$. In case~(ii), $v_1''(0)$ is unrestricted and~$u''(0)$ can take two values.   
\end{proof}

\begin{remark}\label{central1}
These `test' sections are precisely those evaluating $\bmod{\,\tau^2}$ to the image of $1\in H^\vee$ in the 
$\tau=0$ fiber of $\cC_3$. When $N\neq 0$, this is the unique singular point contained in that fiber.
\end{remark}

\begin{remark}
Single-valued sections $z(\tau)$ have the form
\[
z(\tau) = \pm \tau^N \exp\zeta(\tau)
\]
with $\zeta(\tau)$ odd. We can get such $z(\tau)$ to meet the conditions in (i) and (ii), but not in (iii). 
However, we can always choose $v_N(\tau)$ to be even (single-valued in $\tau^2$). 
\end{remark}

\subsection{The case $\mathrm{Pin}^+(2)\subset\mathrm{SO}(3)$.} 
\label{so3workout}
The algebra is now simpler. Every non-trivial 
symplectic representation is a double, and abelianizes both~$\cC_3$ and~$\cC_4$. Since we 
doubly-covered the former $H^\vee$,
the coordinates $z, Z$ are now the square roots of the ones from $\mathrm{SU}(2)$, and so the (now even) 
integer~$N$ is replaced by $N/2$ in formulas. The chiral rings are
\begin{equation}\label{caRSO3}
\begin{split}
\caR_3(H/\{\pm1\};E)^W &= \bC\left[\tau^2, \: u_{N/2},\: v_{N/2} 	\right], \\
\caR_3(\mathrm{SO}(3);E) &= \bC[\tau^2,\:\tau^{-2}u_{N/2}, \: v_{N/2}].
\end{split}
\end{equation} 
\begin{proposition}\label{singsings}
We pass from the first ring to the second by requiring regularity on restricted sections: 
\begin{enumerate}\itemsep0ex
\item When $N=0$, $u_{0} \in O(\tau^2)$; equivalently, $z(0)=\pm1$. 
 This gives $v_0(\tau) =\pm 2 \pmod{\tau^2}$. 
\item When $N>0$, $u_{N/2}(\tau) \in O(\tau^2)$; equivalently, $Val_\tau(z) = 1$ or $N-1$.  
This gives $v_N(\tau) \in O(\tau)$. 
\end{enumerate}
\end{proposition}
\begin{proof} The check proceeds as in the previous proposition.
\end{proof}
\begin{remark}\label{central2}
The test sections~$(u_{N/2},v_{N/2})$ are those evaluating at $\tau=0$ to the center 
of the dual group $\mathrm{SU}(2)$ when $N=0$, and to its image --- the unique singular point 
in the fiber --- for $N>0$. 
\end{remark}

\begin{remark}
The space $\cC_3$ for~$\mathrm{SO}(3)$ can be changed by switching on $w_3\in H^3(BG;\bZ/2)$ (Theorem.~\ref{summaryprop}), 
changing the Weyl involution to $Z\leftrightarrow -Z^{-1}$. This can be absorbed by the change 
$Z\mapsto iZ$, when $\sqrt{-1}$ is incorporated in the coefficient ring.
\end{remark}

\subsection{The case $S_2 \wr \mathrm{U}(1) \subset \mathrm{U}(2)$.} \label{u2workout}
Weights of~$E$ that are \emph{not} multiples of the 
root do not affect the localizations 
of~$\cC_{3,4}$ along the diagonal divisor $\tau: =\tau_1-\tau_2$. Since any $E$ is again polarizable,  
direct calculation (e.g. in~\cite{telc}) shows this case to match the constraints coming from 
$\bP\mathrm{U}(2) =\mathrm{SO}(3)$ under the natural projection. With the modified coordinates 
\[
Z_{1} = \tau^{-N/2}z_1,\quad Z_2 = \tau^{N/2}z_2,
\]
in terms of which
\[
u_{N/2} = \tau^{N/2+1} \left(Z_1 - (-1)^{N/2}Z_2\right),\quad 
	v_{N/2} = \tau^{N/2}\left(Z_1 + (-1)^{N/2} Z_2\right).
\]
the rings \emph{localized along~$\tau=0$} are
\[
\begin{split}
\caR_3(H; E)^{S_2}_\tau = R[u_{N/2}, v_{N/2}]_\tau \\
\caR_3(\mathrm{U}(2); E)_\tau = R[\tau^{-2} u_{N/2}, v_{N/2}]_\tau
\end{split}
\] 
where we set $R: = \bC[\tau_1+\tau_2, \tau_1\tau_2, (z_1z_2)^{\pm1}]$ for convenience. 

\subsection{Root wall correction.} 
We now finally proceed to describe the nonabelian chiral ring. For a root~$\alpha$, with wall $M_{\alpha}$ 
and affine wall $M_\alpha'$, define integers $N_{\alpha}, N_{\alpha}'$ as in~\S\ref {su2workout}: 
$N_\alpha$ is the sum of  coefficients of the multiples of $\alpha$ present in~$E$, and (for long 
roots in type $C$ only) $N_\alpha'$ is the sum of \emph{integral} coefficients. Only the latter are 
relevant to the affine wall~$e^{\alpha/2} = -1$.

The normalizer $G(\alpha)$ and centralizer~$L_\alpha$ of the $\alpha$-wall~$M_\alpha$ 
have one of the forms 
\begin{equation}\label{su2cent}
\begin{matrix}
G(\alpha) & = & \text{(i) }G'\times \mathrm{SU}(2), & \text{(ii)   } G'\times 
	\mathrm{SO}(3), &\text{(iii)    } G'\times_{\bmu_2} 
	\mathrm{SU}(2), \\
L_\alpha &= &\text{(i')  } H'\times \mathrm{SU}(2), & \text{(ii')  } 
	H'\times \mathrm{SO}(3), &\text{(iii')  } 
	H'\times_{\bmu_2}\mathrm{SU}(2)
\end{matrix}
\end{equation}
for a subgroup $G'\subset G$ and the subtorus~$H'\subset H$ whose Lie algebra is orthogonal to $\alpha$ in $\frh$.

The classes of~$\chi$ and~$\kappa$ need not vanish and the rational $H^\vee$-torsors 
over $\frh/W$ and $H/W$ need not be trivial. (For an example, see~\S\ref{notbirat}.) They are 
trivial for $\mathrm{SU}(2), \mathrm{SO}(3)$ and~$\mathrm{U}(2)$, as we saw (for the first, when 
$w_4(E)=0$). We press this to trivialize $\chi_\alpha, \kappa_\alpha$, locally along each (affine) root wall~$M_\alpha$, and affine wall $M_\alpha'$. 
These trivializations replace the factor~$\tau^N$ of the base cases in defining the modified coordinate~$Z$, 
reducing us to the standard presentation of chiral rings. This allows us to impose well-defined vanishing 
conditions for our test sections on the walls.

\begin{lemma}[Sections]\label{sqroots}
Let the group~$G$ be connected. Cancelling the obstructions allows us to choose, locally along each $\alpha$-root 
wall, trivializations $r_\alpha, q_{\alpha}: \frh \to H^\vee$ of the corrected co-cycles $\chi(s_\alpha)$ 
and $\kappa(s_\alpha)$:
\[
\begin{split}
r_{\alpha}(\xi)\cdot s_\alpha r_{\alpha}^{-1}\left(s_{\alpha}\xi\right) = 
	\chi_\alpha(\xi), \\
q_{\alpha}(\xi)\cdot s_\alpha q_{\alpha}^{-1}\left(s_{\alpha}\xi\right) = 
	\kappa_\alpha(\xi). 
\end{split}
\]
Moreover, the projections $r_\alpha^{h_\alpha},q_\alpha^{h_\alpha}$ to the duals of~$H/H'$ are well-defined modulo functions 
for the action of the reflections orthogonal to $s_\alpha$. 
\end{lemma}
\begin{remark}
In correcting the cocycles, we must use  obstruction cancellations coming from~$G$, not from the normalizer $N(H)$. 
An example of faulty cancellation is found in \S\ref{wrongweyl}.\end{remark}

\begin{proof}
Symplectic representations of $L_\alpha$ in (ii') and (iii') are polarizable. There are no ambiguities 
in their Coulomb branches $\cC_{3,4}(L_\alpha;E)$: see Remark~\ref{choices}.  
In case (i), the $H'$-invariant part of $E$ comes from an $\mathrm{SU}(2)$-representation and is 
handled explicitly as before, while the complement is polarizable. So, locally along the root wall, the co-cycles $\chi_\alpha, \kappa_\alpha$ 
are cohomologous to 
\[
\tau^{N_\alpha} = \langle \alpha | \xi\rangle^{N_\alpha\cdot \alpha}, \quad \tau^{N'_\alpha} = \langle \alpha | \xi\rangle^{N'_\alpha\cdot \alpha},
\]
expressible as co-boundaries as in the three base cases discussed. 

Different choices of polarization will give trivializations differing by products of terms the form 
\[
\langle \nu + k\alpha| \xi\rangle^{\nu+k\alpha}\cdot
	\langle \nu - k\alpha| \xi\rangle^{\nu-k\alpha}, 
\]
with $\nu \perp\alpha$ and various half-integral $k$; these equal $\langle \nu |\xi\rangle^{2\nu}$ on 
the root wall, whose $\exp(h_\alpha)$ projection is invisible. The standard action of the Weyl group 
of~$G'$ would change the cocycle to a cohomologous one, and will then preserve the 
$r_\alpha^{h_\alpha}$-projections on the wall. Invariance can then be extended formally, in the normal 
directions to the walls, by averaging over that Weyl group.  
\end{proof}

\subsection{Test sections.}
Let $C_\alpha\subset L_\alpha^\vee$ be the center of the dual group. A \emph{test section $T$} is a 
local section of $\cC^1_{3,4}(L_\alpha; E)$ along $M_{\alpha}$ or $M_\alpha'$ having order of 
contact  $\langle\alpha|\xi\rangle^\alpha$ (respectively, $(1-x^{-\alpha})^\alpha$) with the image of $C_\alpha$. 
Explicitly, for cases (i) and (ii) in~\eqref{su2cent}, we project $L_\alpha^\vee$ to the 
non-abelian factor, whereupon we apply the conditions in 
Propositions~\ref{sings} and~\ref{singsings}. 
In case (iii), $L_\alpha \cong H'' \times \mathrm{U}(2)$ for a suitable $H''\subset H'$ and we 
project to $\mathrm{U}(2)$; from~\S\ref{u2workout}, we may as well project to $\mathrm{SO}(3)$.  
(These two ways to define test sections are reconciled in Remarks~\ref{central1} 
and~\ref{central2}.)

\setcounter{maintheorem}{1}
\begin{maintheorem}\label{algchiral}
The chiral ring~$\caR_{3,4}(G;E)$ is a super-ring of $\caR^1_{3,4}(N(H); E)$, contained within the latter's 
localization away from the walls. It is obtained by including, for each (affine) root $\alpha$,  those 
functions in $\caR^1_{3,4}(N(H); E)[\langle\alpha | \xi\rangle^{-1}]$ which evaluate regularly 
on the $\alpha$-test sections.
\end{maintheorem}

\begin{remark}
When $N_\alpha, N_\alpha' \le 1$, we may choose $s_\alpha$-invariant test sections; otherwise, 
they are only single-valued on $\frh$.
\end{remark}

\begin{proof}
Away from root hyperplanes, this follows from equivariant localization to~$H$. Localization 
to~$H'$ also shows that, generically along an $\alpha$-root wall, the $\cC_{3,4}$ are 
isomorphic under restriction 
to those of the wall normalizer $G(\alpha)$, and the Theorem now follows from the three base cases 
augmented by the Lemma on sections. 
\end{proof}

\section{Abelianization}\label{abel}

When the $H$-restriction of~$E$ contains the doubled representation $\frg_\bH$, we can build the 
Coulomb branches in two steps. Preliminary de-suspension by the index bundles of $N_\frg$ disconnects 
the Bruhat strata of~$\Omega G$, as in Proposition~\ref {splitprop}. If $\frg_\bH$ is a $G$-subrepresentation, this disconnection happens $G$-equivariantly; 
in general, we can only have $N(H)$-equivariance, but this  suffices. 

This gives an (additive) identification of the resulting spectrum with the disjoint union of stabilizers $BL$ of 
one-parameter subgroups. Subsequent de-suspension by $E\ominus \frg_\bH$ leads to an additive equivalence 
advertised in Theorem~\ref{three},
\[
\cC_{3,4}(G;E) = \cC_{3,4}(H;E\ominus\frg_\bH)/W.
\] 
In the homology statement, we must invert the order of the Weyl group to relate $L$-equivariant 
homology with the Weyl invariants in $BH$.

Equality of the multiplicative structures is enforced by localizing away from the root 
hyperplanes on the Toda base and the standard localization theorem to the maximal torus.   

In the more general case, when half-integral multiples of (some of) the roots of~$G$ appear among the weights, 
we can force a partial abelianization of $\cC_3$ by removing corresponding factors from $C_+$ 
and the Weyl-modifying~$\chi$ factors.

\setcounter{section}{0}
\renewcommand{\thesection}{\Alph{section}}
\section{Obstructions for connected groups}
We discuss the obstructions to the construction of chiral rings  for a 
\emph{connected} compact group $G$ with quaternionic representation $E$. The main results are Theorems~\ref{whenodd} and~\ref{kobst}.

\subsection{Quaternionic irreducibles.} \label{quatclass}
Quaternionic representations are self-dual over $\bC$, and an irreducible self-dual complex representation 
is either orthogonal or quaternionic. The simply connected simple groups carrying complex-irreducible 
quaternionic representations are: 
\begin{equation}\label{quatclassif}
\mathrm{Sp},\quad \mathrm{SU}(2\!\!\! \mod{4}), \quad \mathrm{Spin}(\pm 3\!\!\! \mod{8}),\quad 
\mathrm{Spin}(4\!\!\!\mod{8}), \text{ and} \quad E_7.
\end{equation} 
Except for $\mathrm{SU}$, all complex representations of the 
listed groups are self-dual. For $\mathrm{Spin}$ groups, quaternionic are precisely those 
complex-irreducibles that do not factor through 
$\mathrm{SO}$; for the other groups, the 
test is that the unique\footnote{The center of 
$\mathrm{Spin}(4\!\!\!\mod{8})$ is $\bmu_2^{\times 2}$, with the two factors interchanged by the outer 
automorphism; $\mathrm{SO}$ is the quotient by the diagonal $\bmu_2$.} central element of order $2$ should act as $(-1)$.

\begin{proposition}\label{simplegroups}
Let $E$ be a symplectic representation of a simple group $G$.
\begin{enumerate}\itemsep0ex
\item If $w_4(E)\neq0$, then $G$ is a symplectic group.
\item If $E$ is complex-irreducible with $G=\mathrm{Sp}(m)$, then $\dim_\bH E = m\cdot c_2(E)\pmod{4}$.
\item If $E$ is complex-irreducible, $\dim_\bH E$ is odd iff $G=\mathrm{Sp}(m)$ with $m$ and 
$c_2(E)$ both odd.
\end{enumerate}
\end{proposition}

\subsection{First obstruction.} 
Moving to general connected groups, $w_4$ modulo squares is additive on symplectic representations, 
and vanishes for  polarized ones; so it suffices to understand complex-irreducibles. Since any central torus 
acts trivially on those, we may assume that $G$ is semi-simple.

\begin{theorem}\label{whenodd}
For an irreducible symplectic representation $E$ of $G$, 
$w_4(E)=0$ except in one of the 
following (mutually exclusive) cases:
\begin{enumerate}\itemsep0ex
\item $G= G_o\times  \mathrm{Sp}(m) $ and $E = R\otimes S$, with an odd-dimensional orthogonal 
representation~$R$ of $G_o$ and a symplectic representation  
$S$ of $\mathrm{Sp}(m)$ with odd $c_2$; $w_4(E)\neq0$ on $Sp$ and vanishes  on $G_o$.
\item $G= G_o \times_{\bmu_2} \mathrm{Sp}(m) $ and $E = R\otimes S$, with an orthogonal, $(4n+2)$-dimensional 
representation $R$  of $G_o$ and an odd $\bH$-dimensional 
symplectic representation $S$ of $\mathrm{Sp}(m)$.  
\end{enumerate}
In either case, $R$ is orthogonal on each simple factor of $G_o$, so this factorization of~$G$ is unique.
\end{theorem}
\noindent
Note that in case~(ii), $m$ must be odd, and $\bmu_2$ must act via the sign on $R$ and $S$ separately. 
Also, $w_4(E)$ has a square root $\bar{r}\in H^2(BG;\bZ/2)$ with a preferred lift $\bmod{4}$, from the 
larger group $G':=\mathrm{Spin}(4n+2) \times_{\bmu_4} \mathrm{Sp}(m)$ acting on $R\otimes S$.

\begin{corollary}\label{whensqroot}
For a general symplectic $E$, $w_4(E)$ has a square root in~$H^2(BG;\bZ/2)$ iff $E=E'\oplus E''$, 
where the irreducibles in $E'$ admit square roots of $w_4$, while those in $E''$ come in pairs, 
with the irreducibles in each pair factoring as in (i) above with the same $\mathrm{Sp}(m)$. \qed
\end{corollary}

\subsection{Second obstruction.}
The dichotomy in Theorem~\ref{whenodd} shows that a square root~$\bar{r}$ of $w_4$, if it exists, always admits 
a~$\bmod{4}$ lift.
From a trivialization $w_4-\bar{r}^2 = \delta t$, we define the secondary obstruction 
$\sigma:=Sq^2t +\eta\cdot w_4 $ by exploiting the relation $Sq^2w_4 = \delta(\eta\cdot w_4)$ on $BSp$:   
\begin{equation}
\delta (Sq^2t +\eta\cdot w_4) = Sq^2\delta t + \delta(\eta\cdot w_4) = Sq^2\delta t + Sq^2w_4 = (Sq^1\bar{r})^2 =0. 
\end{equation} \label{mod4}
In Theorem~\ref{whenodd}.ii, $\sigma=0$ because its home $H^5/Sq^2H^3$ is zero in~$G'$ (cf.~\S\ref{leraytrans}.h below). 
For the same reason, $\sigma$ vanishes for the summand $E''\subset E$  in Corollary~\ref{whensqroot}: 
$H^5$ comes from~$BG_o$, which in Theorem~\ref{whenodd}.i carries a polarizable (hence unobstructed) 
representation. The Corollary implies that, in investigating the secondary obstruction, 
we may assume that $E$ is irreducible.

\begin{theorem}\label{kobst}
Let $G$ be connected and $E$ an irreducible symplectic representation for which 
$w_4(E)$ has a square root. Then, $\sigma \in H^5(BG;\bZ/2)/Sq^2H^3$ vanishes, except possibly when 
\begin{itemize}\itemsep0ex
\item $G=G_o \times_{\bmu_2} \mathrm{Sp}(m)$, with $m$ odd and a connected group~$G_o$ containing 
a central~$\bmu_2$, \item $E=R\otimes S$, with a $4k$-dimensional orthogonal representation~$R$ of~$G_o$ 
and a quaternionic representation~$S$ of $\mathrm{Sp}(m)$ 
of odd $\bH$-dimension. 
\end{itemize}
In that case: calling~$x$ the generator of~$H^2(B^2\bmu_2; \bZ/2)$ and $\underline{w}_3$ either choice of 
class in $H^3(B\bP\mathrm{SO}(R);\bZ/2)$ lifting to $w_3(R)$ in $SO(R)$,  
\[
\sigma =x\cup \underline{w}_3(R) \pmod{Sq^2 H^3}. 
\]
\end{theorem}
\begin{remark}
The obstruction $\sigma$ and its Bockstein $B\sigma$ are 
non-zero for  $G_o=\mathrm{SO}(4k), E = \bR^{4k}\otimes \bH^m$. The smallest example, 
$\mathrm{SU}(2)^{\times 3}/\mathrm{S}\{\pm1\}^{\times 3}$ 
acting on $(\bC^2)^{\otimes 3}$, is 
discussed in \S\ref{s2not0}.
\end{remark}

\begin{remark}
The map $H^5/Sq^2H^3 \xrightarrow{\ \eta^2\ } KSp^{-1}(BG)$ is well-defined on the $6$-skeleton:  
the next contribution involves~$H^7$. It captures $\eta\otimes E$, there, since the prior cohomological 
term  $\eta\cdot c_2(E)$ vanishes. The image of $Sq^2H^3$ is killed by the Atiyah-Hirzebruch~$d_2$, 
matching the stated ambiguity in~$\sigma$. This collapse also removes the indeterminacy~$x Sq^1x$ 
coming from $\underline{w}_3$: $xSq^1x +Sq^2Sq^1x$ 
is killed by the Leray $d_5\, c_2$ on odd $\mathrm{Sp}(m)$ (\S\ref{leraytrans}.a below).  
\end{remark}

\subsection{Technicalities on $H^4(BG)$.} \label{leraytrans} 
Classes $x\in H^4(BG;\bZ)$ are represented by invariant quadratic forms $q(x)$ on the Lie 
algebra $\frg$ which are integer-valued on the co-weights. We study the Leray spectral sequence for the 
fibration 
\[
B\widetilde{G} \hookrightarrow BG \twoheadrightarrow B^2\pi, 
\]
with $\pi$ finite and $G=\widetilde{G}/\pi$.  We will use the isomorphism 
\[
H^5(B^2\pi;\bZ) \cong  H^4(B^2\pi;\bQ/\bZ);
\] 
by \cite{mac}, the right group classifies $\bQ/\bZ$-valued homogeneous quadratic forms on $\pi$,  
the ratio of the co-weight lattices. 

When $\widetilde{G}$ is simply connected, the leading differential $d_5x\in H^5(B^2\pi;\bZ)$ gives 
the restriction of $q(x) \pmod{\bZ}$ to $\pi$. In general, 
there are prior Leray differentials 
\[
d_2x\in H^2(B^2\pi; H^3(B\widetilde{G};\bZ)), \quad d_3x \in H^3(B^2\pi; H^2(B\widetilde{G};\bZ)),
\] representing the $\bQ/\bZ$-valued bilinear pairing defined from $q(x)$ on 
$\pi\times(\pi_1\widetilde{G})_{tors}$, respectively $\pi\times (\pi_1\widetilde{G})_{free}$. 
Should these two vanish, $d_5x$ is again represented by the (now well-defined) restriction of
$q \pmod{\bZ}$ to $\pi$.

\begin{enumerate}[label=(\alph*)]
\item For $G=\bP\mathrm{Sp}(m)=\mathrm{Sp}(m)/\{\pm 1\}$, with co-root lattice $\langle\pm\mathbf{e}_i\rangle$, 
$q(c_2) = -\sum x_i^2$ and
\[
d_5 (c_2) = -\frac{m}{4}\in \frac{1}{4}\bZ/\bZ =  
	H^5(B^2\pi;\bQ/\bZ).
\]

\item For $G=\bP \mathrm{SU}(n)$, $d_5(c_2)$ is the generator $\frac{1-n}{2n}$ of $H^5(B^2\bmu_n) 
\cong \bZ/(n\cdot \mathrm{gcd}(2,n))$. 

\item In type $D_l$, the co-roots are $\{\pm\mathbf{e}_i\pm\mathbf{e}_j\}_{i<j}$ and $q(p_1)$ is 
the sum-of-squares; the generating class for $\mathrm{Spin}(2l)$ is $p_1/2$, while its center is 
$\bmu_4$ for $l$ odd and $\bmu_2^{\times 2}$ for $l$ even. Generators of~$\pi_1$ are  
$b_\pm:=\left[\pm1/2,\dots, 1/2\right]$, interchanged by the outer automorphism; let also 
$a:=b_+-b_- = [1,0,\dots,0]$. We have
\[
q\left(\frac{p_1}{2}\right): b_\pm \mapsto \frac{l}{8}, \quad a \mapsto \frac{1}{2} \quad\pmod{\bZ}.
\]

\item For $G=\mathrm{SO}(2l)=\mathrm{Spin}(2l)/\langle a \rangle$,  
$p_1$ is the surviving generator. 
\item For $G=\bP\mathrm{SO}(2l)=\mathrm{SO}(2l)/\{\pm 1\}$, $q(p_1)$ sends each generator $b_\pm$ 
to $l/4\pmod{\bZ}$ and pairs it integrally with $a$, so that
\[
d_2p_1=0 \quad\text{and}\quad d_5p_1 = \frac{l}{4}\pmod{\bZ}.
\]
The surviving $H^4$ generators are $p_1$, $2p_1$ and $4p_1$, respectively, for $l = 0,2$ and $\pm1 \pmod{4}$. 

\item For $G=\mathrm{Spin}(4k)/\langle b_+\rangle$, $d_5(p_1/2) = k/4\pmod{\bZ}$; same for $b_-$.

\item In  $G=\bP\mathrm{SO}(4k)$ with the generating classes $u_\pm\in H^2(B^2\pi_1;\bZ/2)$, we have 
\[
H^5(BG;\bZ) = \bZ/4\oplus\bZ/4\oplus\bZ/2;
\] 
the first two summands are generated by 
the Bocksteins $B_4:\bZ/4\to \Sigma \bZ$ of the Pontrjagin squares $\wp(u_\pm)\in H^4(B^2\bmu_2;\bZ/4)$. 
Matching the quadratic form in (c), 
\[
d_5(p_1/2)= k\cdot B_4\wp(u_++u_-) + B(u_+u_-) \in H^5(B^2\pi_1;\bZ). 
\]  
Reducing $B_4\wp(x)$ mod~$2$ gives $x Sq^1 x + Sq^2Sq^1x$. 

\item A recurring example is $G'= \mathrm{Spin}(4k+2)\times_{\bmu_4} \mathrm{Sp}(m)$. 
The successive cohomologies with~$\bZ/2$ coefficients are spanned by classes
\[
H^*(B^2\bmu_4; \bZ/2) = \langle 1\rangle, 0, \langle x\rangle, \langle y\rangle, 
\langle x^2 \rangle, \langle xy, Sq^2y \rangle, \dots  
\]
with $Sq^1x=0$, $y = \frac{1}{2}B_4x \pmod{2}$. Thus $Sq^1(xy)=0$, so $xy$ is the reduction of the 
integral class. That is killed by the differential $d_5(p_1/2)$ from $\mathrm{Spin}$ (Part c above), 
so it vanishes in $H^5(BG';\bZ/2)$. The latter becomes isomorphic to $Sq^2H^3$, because the~$\mathrm{Spin}$ 
and~$\mathrm{Sp}$ factors do not contribute to~$H^5$ with $\bZ/2$ coefficients.   
\end{enumerate}

\begin{proof}[Proof of Proposition \ref{simplegroups}]
For Part (i), note that $Sq^2=0$ on $H^4(B\mathrm{Sp})$, whereas I claim that $Sq^2\neq 0$ on~$H^4$ 
in the other Lie types\footnote{Minding of course the accidental identifications $A_1/C_1, B_2/C_2$.} in \eqref{quatclassif}, thus forcing $w_4(E)$ to vanish. 

If~$\pi_1 G$ has odd order,~$BG$ is equivalent at the prime~$2$ to a simply connected type 
in the list, where the non-vanishing of $Sq^2$ is known. 
The remaining possibility is $G=\mathrm{Spin}(8n+4)/\langle b_+\rangle$ (or $b_-$). Then, 
\ref{leraytrans}.f shows that pull-back from the base $B^2\bmu_2$ induces an isomorphism
\[
H^4(BG;\bZ/2) \cong H^4(B^2\bmu_2;\bZ/2):
\] 
$d_5(p_1/2)\in H^5(B^2\bmu_2;\bZ)$ is a generator, so that $p_1/2$ does not survive in 
the mod~$2$ Leray sequence. Now, $Sq^2H^4(B^2\mu_2) \to H^6 (B^2\bmu_2;\bZ/2)$ is injective, and 
no degree $5$ class is present to give a kernel when mapping to $H^6(BG)$. 

In Part (ii), $(-I)\in \mathrm{Sp}(m)$ maps to $(-I)\in\mathrm{Sp}(E)$, identifying the two 
$B^2\bmu_2$ bases in the fibrations \ref{leraytrans}.a. The transgressions $m\cdot c_2(E)$ and 
$\dim_\bH E$ are thereby equated mod $4$.

Finally, in Part (iii), $G$ has a central element mapping to $(-I)\in \mathrm{Sp}(E)$. As $\dim_\bH(E)$ 
is odd, $c_2(E)$ transgresses to a generator of $H^5(B^2\bmu_2;\bZ)$: in particular, it is odd in 
$H^4(BG)$, and then $G=\mathrm{Sp}(m)$ with $m$ odd, as per 
Parts (i) and (ii).
\end{proof}

\subsection{Proof of the Theorems.} 
From a general quaternionic representation of~$G$, remove the polarizable summands, which do 
not contribute to the obstructions. Now, a finite cover of $G$ splits into simple factors and a 
torus, over which each irreducible summand of the representation is a tensor product of irreducible 
self-dual factors. (In particular, the torus acts trivially.) Given an irreducible summand~$E$, we 
can choose a \emph{simple} factor~$G_s$ with a quaternionic tensor factor $S$ in $E$; the remaining 
factors combined carry an orthogonal representation $R$, and $E=R\otimes S$. Matching this, 
we write $G = G_o\times_{F} G_s$, for a finite central subgroup $F$ of $G_o\times G_s$ which embeds in $G_s$. 

\begin{proof}[Proof of Theorem~\ref{whenodd}]
Here, we assume that $E=E'$ is irreducible. \\
\emph{Case (i)} $F$ has odd order. For $2$-primary questions, we  lift to the product $G_o\times G_s$, where
\begin{equation}\label{c2formula}
c_2(E) = c_2(S) \dim_\bC R - p_1(R)\dim_\bC S.
\end{equation}    
As $\dim_\bC S$ is even, we need both $\dim_\bC R$ and $c_2(S)$ to be odd, and 
Proposition~\ref{simplegroups} settles this case.

\emph{Case (ii)} $F= \bmu_2\times F'$, with $F'$ odd.
Now, $\bmu_2$ acts via the sign on $R$ and $S$: else, we could descend~$E$ to a product group in Case~(i), 
giving the contradiction $\mathrm{Sp}(m)=G_s/\bmu_2$. In particular, $\dim R =2l$, even. Passing 
to $G':= \mathrm{SO}(2l) \times _{\bmu_2} \mathrm{Sp}(S)$, Leray over $B^2\bmu_2$ 
gives an exact sequence 
\begin{equation}
\begin{aligned}
0\to H^4(BG';\bZ)& \to H^4(B\mathrm{SO} \times B\mathrm{Sp};\bZ) = \bZ p_1\oplus \bZ c_2 \xrightarrow{\ d_5\ } \bZ/4 = H^5(B^2\bmu_2;\bZ), \\
d_5&: p_1(R) \mapsto l, \quad  c_2(S)\mapsto (-\dim_\bH S) \pmod{4}.
\end{aligned}
\end{equation}
(As in \S\ref{leraytrans}.e, $d_2=0$: the dot product pairing on $\bmu_2\times 
\pi_1\mathrm{SO}(2l)$ is integral.)  From \eqref{c2formula},  
\begin{equation} \label{c2mixed}
d_5\, \frac{c_2(E)}{2} = d_5\left(c_2(S)\cdot l - p_1(R)\cdot \dim_\bH(S)\right) = -2l\cdot\dim_\bH S \pmod{4}, 
\end{equation}
which is non-zero precisely when both $\dim_\bH S$ and $l$ are odd. In particular, $G_s=\mathrm{Sp}(m)$ 
with $m$ odd, as per Proposition~\ref{simplegroups}.iii, and~$F'$ must be trivial.

\emph{Case (iii)} If $G_s=\mathrm{Spin}(8k+4)$ and $F=\bmu_2^2$, then one of the two $\pm1$ factors 
must act trivially on~$R$, thus also on $S$, and dividing it out in both factors gets a contradiction 
with Cases (i) or (ii). 

In the list \eqref{quatclassif}, $G_s$ never contains a central $\bZ/4$, so we have covered all cases.
\end{proof}

\begin{proof}[Proof of Theorem~\ref{kobst}]
Factor as before $E=R\otimes S$ over $G=G_o \times_F  G_s $. Any central torus in~$G$ must act 
trivially and may be quotiented out, and we may also ignore the odd part of $F$. On $R$ and~$S$, 
$F$ acts either trivially, or else by the same sign. 

\textbf{1.} If $F$ acts trivially, we factor through $F\times F$:
\[
G \xrightarrow{}G_o/F \times  G_s/F \xrightarrow{}\mathrm{SO}(R) \times \mathrm{Sp}(S).
\] 
If $\pi_1(G_s/F)$ is odd, $H^5(BG;\bZ/2)$ comes from the factor $G_o/F$, on which $R\otimes S$ 
is polarizable, and~$\sigma=0$. Else, the 
classification \S\ref{quatclass} ensures that $G_s/F = \mathrm{Spin}(8k+4)/\langle b_\pm\rangle$, 
for one of the generators of the center (\S\ref{leraytrans}.c), with the other generator acting by 
the sign. The~$\mathrm{Spin}$ extra-special~$2$-group ensures that $\dim_\bH S$ is even 
(divisible by $2^{4k}$ in fact). Then, $w_4(R\otimes S)=0$ in the target group, so~$\sigma$ is defined 
and vanishes there by the previous argument.

\textbf{2.} If $F$ acts by a sign, $E$ comes from 
$G':= \mathrm{SO}(R) \times_{\bmu_2} \mathrm{Sp}(S) $, with $\dim R$ even. If $4 \nmid \dim R$, 
$Sq^2H^3= H^5(BG';\bZ/2)$ by \S\ref{leraytrans}.h,  
and again there is no home for $\sigma$. 
The only remaining case is $G=\mathrm{SO}(4k) \times_{\bmu_2} \mathrm{Sp}(m)$, discussed next.

\textbf{3.} 
Denote the generators of the homotopy groups $\bZ/2, \bZ/2, 0,\bZ$ of $KO$ by $\mathbf{1}, \eta, \eta^2$, 
and $\alpha$; subsuming the Bott periodicity generator $\beta\in \pi_8$, $\alpha$ will also denote 
the generator of $\pi_0KSp$, and $\mathbf{1}$ that of $\pi_4BSp$. Recall also from~\S\ref{leraytrans}.g 
the $H^2(B^2\pi_1G;\bZ/2)$ generators $u_+, u_-$, and $x:=u_++u_-$ 
of~$H^2(B^2\bmu_2)$. Either class $u_\pm$ 
lifts to $w_2\in H^2(BSO(4k)$. 

The representations $R, S$ define classes in the $x$-twisted $KO$-groups 
\[
[R]\in {}^x KO^0(\bP\mathrm{SO}(4k)),\quad [S]\in {}^x KSp^0\left(\bP\mathrm{Sp}(m)\right), 
\]
multiplying naturally to the untwisted $KSp^0(BG)$. Denote by ${}^x\alpha$ the (partial) extensions 
of~$\alpha$ to these twisted groups (see the Lemma below). Since~$\dim R$ is even, 
\begin{equation}\label{etaR}
\eta\otimes [R] = B u_\pm\otimes {}^x\alpha \in {}^x KO^{-1}\left(Sk_6 BG\right),
\end{equation}
well-defined, as the ambiguity $Bx \cdot {}^x\alpha$ is killed in  $B\bP\mathrm{SO}$ by $d_3\eta^2$ in 
the twisted Atiyah-Hirzebruch sequence for ${}^xKO$, and the $KO$-theoretic tail is not seen on 
the $6$-skeleton. Similarly, $S$ defines a class which (modulo to the fourth skeletal filtration) is
\[
[S] = {}^x\alpha\cdot \dim_\bH S 
\]
Therefore,
\[
\eta\otimes [R]\otimes [S]= \dim_\bH S \cdot \underline{w}_3(R) \cdot ({}^x\alpha)^2
\]
having denoted by $\underline{w}_3$ either of $Sq^1 u_\pm$. In light of the Lemma that follows, where the 
leading term kills $\eta$, this gives $\dim_\bH S \cdot \underline{W}_3(R)\cup x\cdot \eta^2$, proving the theorem.  
\end{proof}

\begin{lemma} 
For a space $X$ and $x\in H^2(X; \bZ/2)$, denote by ${}^x\alpha$  the unique extension of~$\alpha$ to a 
class in ${}^xKO^4(Sk_3X)$. Then,  
\[
({}^x\alpha)^2 = 4 + \eta^2 x.
\]
\end{lemma}
\begin{proof}
This is a universal formula, so an instance suffices. We use 
\[
X=B\mathrm{SO}(4) = B\left(\mathrm{SU}(2)\times_{\bZ/2} 
	\mathrm{SU}(2)\right)
\]
noting that the tensor product $\bH\otimes_\bH \bH= \bR^4$ of standard representations starts as $4 + w_2\eta^2$.  
\end{proof}

\section{Examples}

\subsection{Obstructed $\cC_4$ but not $\cC_3$.}\label{s2not0}
We describe the smallest connected example exhibiting a secondary obstruction $B\sigma$ for $\cC_4$, 
but not for $\cC_3$ (see~Theorem~\ref{kobst}). Specifically, we show that the Weyl cocycle~$s^{\otimes 2}$ 
of \S\ref{allweyl} is non-zero; this implies that $\sigma\neq0$, and we show that $B\sigma\neq 0$ as well. Let 
\[
G= \mathrm{SU}(2)^{\times 3}/\mathrm{S}(\{\pm I\}^{\times 3}), 
\quad E=\bH^{\otimes 3}.
\] 
Denoting $\mathrm{S}(\{\pm I\}^{\times 3})$ by $\pi$, consider the fibration
\[
B\tilde{H} \rightarrowtail BN(H) \xdbheadrightarrow{\ p\ } 
	BW\times B^2 \pi.
\]
\begin{itemize}\itemsep0ex
\item The Weyl group $W=\bmu_2^{\times 3}$ has generators
$x_{1,2,3}\in H^1_W(\bZ/2)$;
\item the cover 
$\tilde{H}= \mathrm{Spin}(2)^{\times 3}\subset \mathrm{SU}(2)^{\times 3}$ of the maximal torus $H$ carries the Chern classes $\omega_{1,2,3}$;
\item the second base factor has generators $u_{1,2,3}\in H^2\left(B^2 \pi;\bZ/2\right)$, with 
$u_1+u_2+u_3=0$;
\item $p$  is determined by $(x_i,u_i) = (w_1, w_1^2+w_2)$ 
in each of the three projections $N(H) \to \mathrm{O}(2)$.
\item the shortest Leray differentials are 
\begin{equation}\label{xurel}
d_3: \omega_i \mapsto Sq^1u_i+ x_iu_i + x_i^3,
\end{equation}
implementing the relation $Sq^1w_2 + w_1w_2=0$ in each $B\mathrm{O}(2)$. 
\end{itemize}
\noindent
Theorem~\ref{kobst} asserts that $\sigma = Sq^1 u_1\cdot u_2$:  in light of \eqref{xurel}, this is $x_1u_1u_2 + x_1^3u_2$. 
Its transgression $\tau(\sigma) \in H^4_{N(H)}(H;\bZ/2)$ is 
$x_1\cdot \tau(u_1u_2) + x_1^3 \tau(u_2)$, and its leading term in 
 $H^1_W\left(\Lambda^{\otimes 2}/2\right)$, 
coming from $x_1\cdot \tau(u_1u_2)$, represents the $1$-co-cycle which vanishes on all $(x_1,x_2,x_3)\in W$, unless $x_1=(-1)$; 
in which case it takes the value
\begin{equation}\label{concrete}
4\omega_1\otimes \omega_2 + 4\omega_2\otimes \omega _1 \in \Lambda^{\otimes 2} \pmod{2}.
\end{equation}

Let us check it directly on the Weyl group using Theorem~\ref {weylcocycles}.ii. 
The Chern roots of $E$ are $\sum\varepsilon_i \omega_i$, with independent signs $\varepsilon_i$. Choosing $E_+$ defined by 
$\varepsilon_1 = 1$, the cocycle $s^{\otimes 2}$
vanishes on $\bmu_2^{\times 3}$, unless the first entry is 
$(-1)$, in which case we get, modulo $2\Lambda^{\otimes 2}$, 
\[
\begin{split}
(\omega_1&+\omega_2+\omega_3)^{\otimes 2} + 
(\omega_1+\omega_2-\omega_3)^{\otimes 2} + 
(\omega_1-\omega_2+\omega_3)^{\otimes 2} + 
(\omega_1-\omega_2-\omega_3)^{\otimes 2}  =\\
&= 4\sum\omega_i^{\otimes 2} = 4\omega_1^{\otimes 2} + 4\omega_2^{\otimes 2} + 4(\omega_1+\omega_2)^{\otimes 2} = 4\left(\omega_1\otimes \omega_2 + \omega_2\otimes \omega _1\right) + 
	8\left(\omega_1^{\otimes 2} + \omega_2^{\otimes 2}\right).
\end{split}
\]
The second term vanishes, because $2\omega_i\in \Lambda$; so this matches \eqref{concrete} mod~$2$.  

Then, $\sigma$ must be the unique non-zero class in $H^5(BG;\bZ/2)/Sq^2H^3$, represented by 
$u_iSq^1u_j$ for any pair $i\neq j$; and its Bockstein does not vanish. 

\subsection{A disconnected group with $B\sigma\neq0$.} 
\label{noK}
Abandoning connectivity makes second obstructions easier to find. Let $G$ be the extension of 
$\bmu_2$ by $H:=\mathrm{U}(1)^{\times 3}_-$, with $\bmu_2$ inverting each factor and extended 
via the tri-diagonal class in $H^2(B\bmu_2;H)$. Equivalently, the 
non-trivial component of $G$ squares 
to $(-1)^{\times 3}\in H$.  This is the subgroup of $(\mathrm{Pin}^-(2))^{\times 3}$ lying over the 
diagonal $\bmu_2$; it also carries an extra homomorphism to $\mathrm{Pin}^-(2)$, by multiplying out 
the three $\mathrm{U}(1)$ factors. Consider now the representation 
\[
E = (L_1\oplus L_2 \oplus L_3 \oplus L_1L_2L_3) \oplus (\text{dual})
\]
built from the standard representations $L_i$ of the three factors of $H$, and extended to 
a symplectic representation of $G$ via the four maps to $\mathrm{Pin}^-(2)$. 

Let $\Lambda = H^2(BH;\bZ)$ with the basis of the three Chern classes $\omega_i$, call $[\omega_i]$ 
their reductions mod~$2$, and $u$ the generator of $H^1(B\bmu_2;\bZ/2)$. I claim that
\begin{enumerate}
\item $H^4(BG;\bZ)$ injects into $H^4(BH;\bZ)$;
\item $h:=\omega_1\omega_2+\omega_1\omega_3+\omega_2\omega_3$ is in the image. (More precisely, 
$H^4(BG)=\langle h,\omega_i^2\rangle$.)
\end{enumerate}
Both claims are seen from the Leray sequence for $BH \rightarrowtail BG \twoheadrightarrow B\bmu_2$. 
With $\bZ/2$ coefficients, $d_3[\omega_i] = u^3$, from the central extension. Item (i) holds because $H^2(B\bmu_2; \Lambda)=0$, while 
$H^4(B\bmu_2;\bZ)$ is killed by $d_3(B[\omega_i]) = B(u^3)$. Part (ii) holds because $h$ survives in the integral Leray sequence: reducing mod~$2$ shows that
\[
d_3(\omega_i\omega_j) = u^2\cdot B([\omega_i]+[\omega_j]),
\] 
so that $d_3h=0$, and there is no landing place for $d_5$. We conclude that
\[
c_2(E) = \omega_1^2 + \omega_2^2 +\omega_3^2 + \left(\omega_1+\omega_2+\omega_3\right)^2 = 
	2\left(\omega_1\omega_2+\omega_1\omega_3+\omega_2\omega_3\right) + 2\sum\omega_i^2  =2\left(h + \sum \omega_i^2\right)
\]
and is even in $H^4(BG;\bZ)$. 

Now, we have $Sq^1[\omega_i] = u\cdot[\omega_i]$ and $Sq^2[\omega_i] = [\omega_i]^2 \pmod{u^4}$; in particular, this holds over the $2$-skeleton $\bR\bP^2$ 
of the base. We compute
\[
\begin{split}
&Sq^2 h = Sq^2(\omega_1\omega_2+\omega_1\omega_3+\omega_2\omega_3) = \sum_{i\neq j}\: [\omega_i]^2\cdot [\omega_j] 
+ u^2h \\
&Sq^3h = Sq^1Sq^2h = u\cdot \sum_{i\neq j}\: [\omega_i]^2\cdot [\omega_j] \\
&Sq^3(\omega_i^2) = 0
\end{split}
\]
Since $Sq^2\sigma = Sq^3h$, we conclude that $\sigma = u\cdot h + u\cdot\varphi$, 
with $Sq^2(u\varphi) = 0$ and $\varphi$ quadratic in the $[\omega_i]$, because $u^3=0$. Symmetry allows only 
$0$ or $\sum\omega_i^2$ as options. Then,
\[
B\sigma = u^2\cdot (h+\varphi) \in H^6(BG;\bZ).
\]
Moreover, $H^3(BG;\bZ/2)$ is spanned by the $u\cdot[\omega_i+\omega_j]$, and the class above is not 
a sum of their $Sq^2$, so the obstruction class $B\sigma$ does not vanish modulo the ambiguity $BSq^2 H^3$.

\begin{remark}
We can determine $\sigma$ by restricting to the diagonal $\mathrm{Pin}^-(2) \subset G$. The representation becomes 
\[
L^{\oplus 3} \oplus L^{\otimes 3} \oplus \text{(dual)}
\]
which is restricted from the representation $\bC^2\oplus \bC^2\oplus \bC^4$ of $\mathrm{SU}(2)$, 
and has no obstruction. 
Since $h$ restricts to $3\omega^2$, we need $\varphi = \sum\omega_i^2 = 3\omega^2$ to cancel $B\sigma$. Thus, $\sigma = uh+u\sum[\omega_i^2]$.
\end{remark}

\subsection{A non-trivial Toda torsor.} \label{notbirat}
When the Weyl cocycle $\chi$ is not exact, the space $\cC_3$ is a non-trivial rational torsor over
the Toda group scheme. The smallest example for a connected group may be $G=\mathrm{Sp}(2)$ with the irreducible 
component $E\subset \bR^{5}\otimes \bH^{2}$ complementary to the standard representation $\bH^{2}$. 

Recall that the co-roots are $\pm \mathbf{e}_i$ and $\pm \mathbf{e}_i \pm\mathbf{e}_j$,~$i\neq j$, 
and the weight lattice is spanned by the dual basis $\{\xi_{1,2}\}$ of $\{\mathbf{e}_{1,2}\}$.
Under the normalizer $N(H)$, $E$ splits as two copies of~$\bH^{2}$, plus the sum of the eight 
weight spaces for 
\[
\nu_{\pm 2, \pm 1} = \pm 2\xi_1 \pm \xi_2\quad\text{and}\quad 
\nu_{\pm 1,\pm 2} = \pm \xi_1 \pm 2\xi_2,
\] 
with independent sign choices. Polarize by one copy of $\bH^2$, together with $\nu_{2,\pm1}$ and 
$\nu_{1,\pm2}$ and focus on the inversion 
$(\xi\leftrightarrow -\xi, h\leftrightarrow h^{-1})$. The relevant $\chi$ is 
\[
\chi(\xi_1,\xi_2) = 
	\left[
\begin{matrix}\nu_{1,2}\cdot \nu_{2,1}^2\cdot \nu_{1,-2} \cdot \nu_{2,-1}^2\\
\nu_{1,2}^2\cdot \nu_{2,1}\cdot \nu_{1,-2}^{-2} \cdot \nu_{2,-1}^{-1} 
\end{matrix} 
\right] \in \left[\begin{matrix}
\bC^\times \\ \bC^\times\end{matrix}\right] = H^\vee 
\] 
up to correction by a constant $2$-torsion point in $H^\vee$. This $\chi$ cannot be expressed as 
\[
\chi = \delta\varphi := \varphi(-\xi)^{-1}/\varphi(\xi),
\] 
for any rational map $\varphi:\frh\to H^\vee$: such a~$\delta\varphi$ must have even valuation 
at all of the lines~$\nu_{i,j}=0$, but that is not the case on any of those lines. The torsor~$\cC_3(G;E)$ 
is therefore rationally non-trivial. 

\subsection{Erroneous removal of the Weyl obstruction.} 
\label{wrongweyl}
Consider the group $G=\mathrm{SU}(2)\times \mathrm{U}(1) $, with 
representation $E= \bC^2\otimes \left(\bC\oplus \overline{\bC}\right)$. 
The natural polarization supplies two Lagrangian sections of $\cC_3(G;E)$, giving birational equivalences to the Toda space. 

Let us instead choose $E_+:=\langle \mathbf{e}_1, \overline{\mathbf{e}}_2\rangle$,  
the highest-weight space for $\mathrm{SU}(2)$. Up to possible adjustment by a $2$-torsion point in 
$H= \bC^\times\times \bC^\times$, the Weyl reflection (acting on the first coordinates of $\frh$ and $H$) 
now gets multiplied by the factor 
\begin{equation}\label{badweyl}
\chi(\xi_1,\xi_2) = \left[\begin{matrix}
(\xi_1+\xi_2)(\xi_1-\xi_2)\\ 
(\xi_1+\xi_2)(\xi_1-\xi_2)^{-1}
\end{matrix}\right]
\in \left[\begin{matrix}\bC^\times \\ \bC^\times\end{matrix}\right] \subset \left[\begin{matrix}\mathrm{SL}_2 \\ 
	\mathrm{GL}_1\end{matrix} \right]  
\end{equation}
This modified reflection squares to $1$, so it seems that no correction is needed. However, 
after restricting to the root hyperplane $\xi_1=0$, the natural Weyl action gets multiplied by 
$[-\xi_2^2, -1]$. This gives the non-trivial class in 
$H^1_{\bZ/2}(\bC^\times)$ on the 
second $\bC^\times$ factor, and can thus not be used to define the trivial torsor $\cC_3(G;E)$. 

To apply the Weyl descent construction, we must change the sign in the bottom entry of \eqref{badweyl}. 
This exploits the ambiguity $\bZ/2$ in the construction of $\cC_3(N(H); E)$, the integral Bockstein image of 
$H^3_{N(H)}(\bZ/2)$. By contrast, $\cC_3(G;E)$ is unambiguous: see Remark~\ref{obsremark}.i.

\bigskip
\bigskip
\noindent\small
{CONSTANTIN TELEMAN\\
Dept.~of Mathematics, 970 Evans Hall, UC Berkeley, Berkeley, CA 94720, USA\\ 
\texttt{teleman@berkeley.edu}
}

\end{document}